\documentclass[11pt]{article}

\title{The Weinstein conjecture for stable Hamiltonian structures}
 \author{Michael Hutchings\footnote{Partially supported by NSF grant
     DMS-0505884.} and Clifford Henry Taubes\footnote{Partially
     supported by the National Science Foundation.}}
\date{}

\usepackage{amssymb}
\usepackage{latexsym}
\usepackage{amsmath}
\usepackage{amsthm}
\usepackage{amscd}

\newcommand{\mc}[1]{{\mathcal #1}}

\numberwithin{equation}{section}

\newtheorem{theorem}{Theorem}[section]
\newtheorem{proposition}[theorem]{Proposition}

\newtheorem{lemma}[theorem]{Lemma}
\newtheorem{lemma-definition}[theorem]{Lemma-Definition}

\theoremstyle{definition}
\newtheorem{definition}[theorem]{Definition}

\newtheorem{example}[theorem]{Example}

\newtheorem{notation}[theorem]{Notation}

\newtheorem*{examplestar}{Example}

\newcommand{\floor}[1]{\left\lfloor #1 \right\rfloor}
\newcommand{\ceil}[1]{\left\lceil #1 \right\rceil}

\newcommand{\eqdef}{\;{:=}\;}

\renewcommand{\frak}{\mathfrak}

\newcommand{\C}{{\mathbb C}}

\newcommand{\R}{{\mathbb R}}

\newcommand{\Z}{{\mathbb Z}}
\newcommand{\op}{\operatorname}

\newcommand{\Ker}{\op{Ker}}

\newcommand{\bpm}{\begin{pmatrix}}
\newcommand{\epm}{\end{pmatrix}}

\renewcommand{\epsilon}{\varepsilon}

\begin{document}

\setcounter{tocdepth}{2}

\maketitle

\begin{abstract}
  We use the equivalence between embedded contact homology and
  Seiberg-Witten Floer homology to obtain the following improvements on the
  Weinstein conjecture.  Let $Y$ be a closed oriented
  connected 3-manifold with a stable Hamiltonian structure, and let
  $R$ denote the associated Reeb vector field on $Y$.  We prove that
  if $Y$ is not a $T^2$-bundle over $S^1$, then $R$ has a closed
  orbit.  Along the way we prove that if $Y$ is a closed oriented
  connected 3-manifold with a contact form such that all Reeb orbits
  are nondegenerate and elliptic, then $Y$ is a lens space.  Related
  arguments show that if $Y$ is a closed oriented 3-manifold with
  a contact form such that all Reeb orbits are nondegenerate, and if
  $Y$ is not a lens space, then there exist at least three distinct
  embedded Reeb orbits.
\end{abstract}

\tableofcontents

\section{Introduction}

Throughout this paper, $Y$ denotes a closed, oriented, connected
3-manifold.  Recall that a {\em contact form\/} on $Y$ is a 1-form
$\lambda$ on $Y$ such that $\lambda\wedge d\lambda >0$.  A contact
form $\lambda$ determines a contact structure, namely the two-plane
field $\xi\eqdef \Ker(\lambda)$.  It also determines a vector field
$R$, called the {\em Reeb vector field\/}, characterized by
$d\lambda(R,\cdot)=0$ and $\lambda(R)=1$.  A {\em Reeb orbit\/} is a
closed orbit of $R$, ie a map $\gamma:\R/T\Z\to Y$ for some $T>0$ such
that $\gamma'(t)=R(\gamma(t))$.  Two Reeb orbits are considered
equivalent if they differ by precomposition with a translation of
$\R/T\Z$.

The three-dimensional version of the {\em Weinstein conjecture\/}
asserts that for every closed oriented 3-manifold $Y$, and for every
contact form $\lambda$ on $Y$, there exists a Reeb orbit.  There is a
long history of work proving this conjecture in many cases, for
example for overtwisted contact structures in \cite{h}, for contact
structures supported by planar open books in \cite{ach}, and for many
additional types of open books in \cite{ch}.  The Weinstein conjecture
was recently proved by the second author in all cases in \cite{tw}.

In fact a stronger result is proved in \cite{echswf}, which asserts
that a version of the Seiberg-Witten Floer homology of $-Y$ as defined
by Kronheimer-Mrowka \cite{km}, namely the version $\check{HM}(-Y)$,
is isomorphic to the embedded contact homology (ECH) of $(Y,\lambda)$.
To see how this implies the Weinstein conjecture, recall that the ECH
of $(Y,\lambda)$ is the homology of a chain complex which is generated
by certain unions of Reeb orbits with multiplicities, and whose
differential counts certain embedded holomorphic curves in $\R\times
Y$.  If $Y$ has a contact form with no Reeb orbit, then the ECH has
just one generator, namely the empty set of Reeb orbits.  However
Kronheimer and Mrowka have shown that $\check{HM}$ of every closed
oriented 3-manifold is infinitely generated \cite{km}.  Together with
the isomorphism between ECH and $\check{HM}$, this gives a
contradiction, and the Weinstein conjecture is proved.

The original proof of the Weinstein conjecture in \cite{tw}
established and used only a first step towards the isomorphism between
ECH and $\check{HM}$, which entailed showing how generators of
$\check{HM}$ give rise to generators of ECH.  In the present paper, we
exploit the full strength of the isomorphism to obtain two
improvements on the Weinstein conjecture.  First, we extend the class
of vector fields for which a closed orbit is known to exist.  Second,
in most cases we can increase the lower bound on the number of
(embedded) Reeb orbits from one to three.

To explain the first improvement: A contact form is a special case of
a {\em stable Hamiltonian structure\/}, a notion which was identified
in \cite{behwz,cm} as a general setting in which one can obtain
Gromov-type compactness for moduli spaces of holomorphic curves in
$\R\times Y$.  If $Y$ is an oriented three-manifold, a stable
Hamiltonian structure on $Y$ is a pair $(\lambda,\omega)$, where
$\lambda$ is a $1$-form on $Y$ and $\omega$ is a $2$-form on $Y$, such
that
\begin{gather*}
d\omega  = 0,\\
\lambda\wedge\omega > 0,\\
\Ker(\omega) \subset \Ker(d\lambda).
\end{gather*}
Note that the second condition implies that $\omega$ is nonvanishing,
and consequently the third condition is equivalent to
\[
d\lambda = f \omega
\]
where $f:Y\to\R$ is a smooth function.

A stable Hamiltonian structure determines a vector field $R$ on $Y$,
which we call the {\em Reeb vector field\/}, characterized by
\[
\omega(R,\cdot) = 0, \quad\quad \lambda(R)=1.
\]
It follows from the definitions that the stable Hamiltonian structure
is invariant under $R$, ie
\[
\mc{L}_R\lambda=0,\quad\quad\mc{L}_R\omega=0, \quad\quad \mc{L}_Rf=0.
\]

\begin{examplestar}
  A contact $1$-form $\lambda$ determines a stable Hamiltonian
  structure in which $\omega=d\lambda$ and $R$ is
  the Reeb vector field in the usual sense.  Here $f\equiv 1$.
\end{examplestar}

\begin{examplestar}
  Let $\Sigma$ be a closed oriented surface with a symplectic form
  $\omega$, and let $\phi$ be a symplectomorphism of
  $(\Sigma,\omega)$.  Let $Y$ be the mapping torus
\begin{equation}
\label{eqn:mappingTorus}
Y\eqdef \frac{[0,1]\times\Sigma}{(1,x)\sim(0,\phi(x))}.
\end{equation}
This fibers over $S^1=\R/\Z$.  Let $R$ denote the vector field on $Y$
which increases the $[0,1]$ coordinate in \eqref{eqn:mappingTorus}.
Note that closed orbits of $R$ correspond to periodic orbits of
$\phi$.  Now $Y$ has a stable Hamiltonian structure in which $R$ is as
described above, $\lambda$ is the pullback of the volume form on
$S^1$, and $\omega$ is the unique extension of the closed $2$-form on
the fibers of $Y\to S^1$ to a $2$-form on $Y$ which annihilates $R$.
Here $f\equiv 0$.
\end{examplestar}

The main result of this paper is the following:

\begin{theorem}
\label{thm:main}
Let $Y$ be a closed oriented connected 3-manifold with a stable
Hamiltonian structure.  If $Y$ is not a $T^2$-bundle over $S^1$, then
the associated Reeb vector field has a closed orbit.
\end{theorem}

Note that there exist $T^2$-bundles over $S^1$ with stable Hamiltonian
structures having no Reeb orbit.  Examples with $f=0$ are provided by
the mapping tori of symplectomorphisms $\phi:T^2\to T^2$ with no
periodic orbit.  Here $\phi$ can be an irrational rotation, or an
appropriate composition of an irrational rotation with a Dehn twist.
In addition, the construction in \S\ref{sec:conclusion} can be
inverted to produce examples of $T^2$-bundles over $S^1$ with stable
Hamiltonian structures in which $f$ changes sign and there is no Reeb
orbit.

The outline of the proof of Theorem~\ref{thm:main} is as follows.
Assume that $R$ has no closed orbit.  By the Weinstein conjecture, $f$
is sometimes zero; and if $f$ is identically zero then $Y$ fibers over
$S^1$ and a calculation using the Lefschetz fixed point theorem shows
that the fiber is $T^2$.  So assume that $f$ is sometimes zero and
sometimes nonzero. Choose $\epsilon>0$ small such that both $\epsilon$
and $-\epsilon$ are regular values of $f$.  Since $f$ is invariant
under the nonvanishing vector field $R$, the level sets
$f^{-1}(\pm\epsilon)$ are disjoint unions of tori.  We can now
decompose $Y$ along a union of tori as
\[
Y = f^{-1}(-\infty,-\epsilon] \cup f^{-1}[-\epsilon,\epsilon] \cup
f^{-1}[\epsilon,\infty).
\]
It is enough to show that each piece in the decomposition is a
disjoint union of copies of $T^2\times I$.

To do this for the middle piece $f^{-1}[-\epsilon,\epsilon]$, the idea
is to use the fact that $f$ is close to zero to show that the middle
piece fibers over $S^1$.  We then use the assumption that there are no
Reeb orbits to show that the fibers are disjoint unions of annuli.
This part of the argument does not use any Floer homology.

The remaining pieces in the decomposition can be collapsed into closed
contact manifolds with each boundary torus becoming an elliptic Reeb
orbit.  To conclude the proof of Theorem~\ref{thm:main}, we then
invoke the following theorem.  To state it, recall that a Reeb orbit
$\gamma$ is called {\em nondegenerate\/} if its {\em linearized return
  map\/}, ie the symplectic linear map from $\xi_{\gamma(0)}$ to
itself given by the linearized Reeb flow along $\gamma$, does not have
$1$ as an eigenvalue.  We can classify the nondegenerate Reeb orbits
into three types according to the eigenvalues $\lambda,\lambda^{-1}$
of the linearized return map:
\begin{itemize}
\item
{\em elliptic\/}: $\lambda,\lambda^{-1}$ on the unit circle.
\item
{\em positive hyperbolic\/}: $\lambda,\lambda^{-1}>0$.
\item
{\em negative hyperbolic\/}: $\lambda,\lambda^{-1}<0$.
\end{itemize}

\begin{theorem}
\label{thm:lens}
Let $Y$ be a closed oriented connected 3-manifold with a contact form
$\lambda$ such that all Reeb orbits are nondegenerate and elliptic.
Then $Y$ is a lens space, there are exactly two embedded Reeb orbits,
and they are the core circles in the solid tori of a genus one
Heegaard splitting of $Y$.
\end{theorem}

The proof of Theorem~\ref{thm:lens} extensively uses the isomorphism
between ECH and $\check{HM}$.  The outline is as follows.  If all Reeb
orbits are elliptic, then the differential on the ECH chain complex
vanishes, because all generators have even grading.  The isomorphism
with $\check{HM}$ then implies that all Reeb orbits represent torsion
homology classes. If the number of embedded Reeb orbits is less than
or greater than two, then the growth rate of the number of
nullhomologous ECH generators with grading $\le I$ as $I$ goes to
infinity is either too slow or too fast to be consistent with known
properties of $\check{HM}$.  Next, the equivalence of the $U$ maps in
the two theories guarantees the existence of many holomorphic curves
between the ECH generators.  Finally, a probabilistic argument shows
that some of these holomorphic curves are in fact cylinders whose
projections to $Y$ are embeddings.  These give rise to a foliation of
$Y$ by holomorphic cylinders with boundary on the two Reeb orbits, and
this foliation yields the desired Heegaard splitting.

A corollary of Theorem~\ref{thm:lens} is that it is impossible for
there to exist only one embedded Reeb orbit and for that orbit to be
nondegenerate and elliptic.  Using ECH$=\check{HM}$ a bit more, one can
upgrade this to show:

\begin{theorem}
\label{thm:three}
Let $Y$ be a closed oriented 3-manifold with a contact form such that
all Reeb orbits are nondegenerate.  Then there are at least two
distinct embedded Reeb orbits; and if $Y$ is not a lens space then
there are at least three distinct embedded Reeb orbits.
\end{theorem}

There is certainly room for improvement on the lower bound in
Theorem~\ref{thm:three}, at least if one knows more about the
three-manifold $Y$ and its contact structure.  In fact, work of Colin
and Honda \cite{ch} using linearized contact homology shows that many
three-manifolds with contact structures have the property that for any
contact form there must be infinitely many distinct embedded Reeb
orbits.

\paragraph{The rest of the paper.} In \S\ref{sec:ECHReview} we review
the basics of ECH.  In \S\ref{sec:subtleties} we discuss some more
subtle aspects of ECH that we will need.  In \S\ref{sec:lens} we prove
Theorem~\ref{thm:lens} regarding contact manifolds with all Reeb
orbits elliptic, and we also prove the lower bound in
Theorem~\ref{thm:three}.  In \S\ref{sec:conclusion} we prove the main
Theorem~\ref{thm:main}.

\section{Review of embedded contact homology}
\label{sec:ECHReview}

We now review the basic notions from embedded contact homology (ECH)
that will be needed in the present paper.  References are
\cite{pfh2,ir} for the ECH index theory, \cite{t3} for additional
structure on ECH, and \cite{obg1} for the analysis.  Below, fix a
closed oriented connected 3-manifold $Y$ as usual, and fix a contact
form $\lambda$ on $Y$ such that all Reeb orbits are nondegenerate.

\subsection{ECH generators}

\begin{definition}
  An {\em ECH generator\/} is a finite set of pairs
  $\alpha=\{(\alpha_i,m_i)\}$ such that the $\alpha_i$'s are distinct
  embedded Reeb orbits, the $m_i$'s are positive integers, and $m_i=1$
  whenever $\alpha_i$ is hyperbolic.  The homology class of $\alpha$
  is defined to be
\[
[\alpha]\eqdef \sum_im_i[\alpha_i]\in H_1(Y).
\]
\end{definition}

For each $\Gamma\in H_1(Y)$, the embedded contact homology
$ECH_*(Y,\lambda,\Gamma)$ is the homology of a chain complex
$C_*(Y,\lambda,\Gamma)$ which is freely generated over $\Z$ by the ECH
generators $\alpha$ with $[\alpha]=\Gamma$.

Before defining the grading and the differential on this chain
complex, we make two remarks.  First, the empty set $\alpha=\phi$ is a
legitimate ECH generator with $[\alpha]=0$.  In fact the ECH generator
$\phi$ turns out to be a cycle in the ECH chain complex
$C_*(Y,\lambda,0)$, whose homology class in ECH conjecturally agrees
with the contact invariant in $\check{HM}$.  Second, we sometimes
write an ECH generator $\alpha=\{(\alpha_i,m_i)\}$ using the
multiplicative notation $\prod_i\alpha_i^{m_i}$.  However the grading
and differential on the ECH chain complex do not behave simply with
respect to this sort of multiplication.

\subsection{The ECH index}

We now explain the grading on the chain complex.

\begin{notation}
If $\alpha=\{(\alpha_i,m_i)\}$ and $\beta=\{(\beta_j,n_j)\}$ are two
ECH generators with $[\alpha]=[\beta]=\Gamma$, define
$H_2(Y,\alpha,\beta)$ to be the set of equivalence classes of
$2$-chains $Z$ in $Y$ with
\[
\partial Z =
\sum_im_i\alpha_i-\sum_jn_j\beta_j,
\]
where two such 2-chains are considered equivalent if they differ by
the boundary of a 3-chain.  Thus $H_2(Y,\alpha,\beta)$ is an affine
space over $H_2(Y)$.
\end{notation}

\begin{definition}
If $Z\in H_2(Y,\alpha,\beta)$, define the {\em ECH index\/}
\[
I(\alpha,\beta,Z) \eqdef c_\tau(Z) + Q_\tau(Z) +
\sum_i\sum_{k=1}^{m_i}\op{CZ}_\tau(\alpha_i^{k}) -
\sum_j\sum_{k=1}^{n_j}\op{CZ}_\tau(\beta_j^k) \in \Z.
\]
Here $\tau$ is a trivialization of the contact plane field $\xi$ over
the $\alpha_i$'s and $\beta_j$'s; $c_\tau(Z)$ denotes the relative
first Chern class over $Z$ with respect to the boundary trivialization
$\tau$; $Q_\tau(Z)$ denotes the relative intersection pairing; and
$\op{CZ}_\tau(\gamma^k)$ denotes the Conley-Zehnder index of the
$k^{th}$ iterate of $\gamma$.  For the detailed definitions of the
integers $c_\tau(Z)$, $Q_\tau(Z)$, and $\op{CZ}_\tau(\gamma^k)$ see
\cite{pfh2,ir}.
\end{definition}

As explained in \cite{pfh2,ir}, the ECH index $I(\alpha,\beta,Z)$ does
not depend on the choice of trivialization $\tau$, even though the
individual terms in its definition do.  It does depend on $Z$:
if $Z'\in H_2(Y,\alpha,\beta)$ is another relative homology clas, then
by \cite[Eq.\ (5)]{pfh2} and \cite[Lem.\ 2.5(a)]{pfh2} we have
\begin{align}
\label{eqn:cZZ'}
c_\tau(Z) - c_\tau(Z') &= \langle c_1(\xi),Z- Z'\rangle,\\
\label{eqn:QZZ'}
Q_\tau(Z) - Q_\tau(Z') &= 2\Gamma \cdot (Z-Z').
\end{align}
Therefore
\[
I(\alpha,\beta,Z) - I(\alpha,\beta,Z') = \langle c_1(\xi) + 2
\op{PD}(\Gamma),Z-Z'\rangle.
\]
Consequently, if $[\alpha]=[\beta]=\Gamma$, then it makes sense to
define
\begin{equation}
\label{eqn:Ialphabeta}
I(\alpha,\beta) \eqdef I(\alpha,\beta,Z)\in \Z/d_\Gamma,
\end{equation}
where $Z$ is any element of $H_2(Y,\alpha,\beta)$, and $d_\Gamma$
denotes the divisibility of $c_1(\xi) + 2\op{PD}(\Gamma)$ in
$H^2(Y;\Z)$ mod torsion.  Note that $d_\Gamma$ is an even integer.

It is also shown in \cite[Prop.\ 1.6(b)]{pfh2} that $I$ is additive in
the following sense: If $\gamma$ is another ECH generator with
$[\gamma]=\Gamma$ and if $W\in H_2(Y,\beta,\gamma)$, then $Z+W\in
H_2(Y,\alpha,\gamma)$ is defined and
\[
I(\alpha,\gamma,Z+W) = I(\alpha,\beta,Z) + I(\beta,\gamma,W).
\]
It follows that \eqref{eqn:Ialphabeta} defines a relative $\Z/d_\Gamma$
grading on the chain complex $C_*(Y,\lambda,\Gamma)$.

It is further shown in \cite{ir} that the relative grading
\eqref{eqn:Ialphabeta} can be refined to an absolute grading which
associates to each ECH generator a homotopy class of oriented 2-plane
fields on $Y$.  In the present paper we will not need this absolute
grading and can just regard the grading on $C_*(Y,\lambda,\Gamma)$ as
taking values in some abstract affine space over $\Z/d_\Gamma$.
However we do need to know, from \cite[Prop.\ 1.6(c)]{pfh2}, that the
mod 2 grading is given by
\begin{equation}
\label{eqn:I2}
I(\alpha,\beta) \equiv I_2(\alpha) - I_2(\beta) \in \Z/2.
\end{equation}
Here if $\alpha=\{(\alpha_i,m_i)\}$ is an ECH generator, then
$I_2(\alpha)\in\Z/2$ denotes the mod 2 count of orbits $\alpha_i$ that
are positive hyperbolic.

\subsection{The index inequality}

To prepare to define the differential on the ECH chain complex, choose
an almost complex structure $J$ on $\R\times Y$ satisfying the
following properties.  Let $s$ denote the $\R$ coordinate on $\R\times
Y$, and recall that $R$ denotes the Reeb vector field.  We require
that $J(\partial/\partial s) = R$, that $J$ is invariant under the map
$(s,y)\mapsto (s+c,y)$ for each $c\in\R$, and that $J$ sends the
contact plane field $\xi$ to itself, rotating positively in the sense
that $d\lambda(v,Jv)\ge 0$ for all $v\in\xi$.

We consider holomorphic curves $u:(C,j)\to(\R\times Y,J)$ such that
the domain $C$ is a punctured compact Riemann surface, and $u$ is not
constant on any component of $C$. The domain $C$ is not required to be
connected. We often abuse notation and denote the holomorphic curve
$u$ simply by $C$.

If $C$ is a holomorphic curve in $\R\times Y$, if $\gamma$ is an
embedded Reeb orbit, and if $k$ is a positive integer, then a
``positive end of $C$ at $\gamma$ of multiplicity $k$'' is an end of
$C$ which is asymptotic to $\R$ cross the $k$-fold iterate of $\gamma$
as $s\to\infty$.  A ``negative end'' is defined analogously but with
$s\to-\infty$.

If $\alpha=\{(\alpha_i,m_i)\}$ and $\beta=\{(\beta_j,n_j)\}$ are two
ECH generators, let $\mc{M}(\alpha,\beta)$ denote the moduli space of
holomorphic curves in $\R\times Y$ with positive ends at $\alpha_i$ of
total multiplicity $m_i$, negative ends at $\beta_j$ of total
multiplicity $n_j$, and no other ends.  If $C\in\mc{M}(\alpha,\beta)$,
then the projection of $C$ to $Y$ has a well-defined relative homology
class $[C]\in H_2(Y,\alpha,\beta)$.  We write $I(C)\eqdef
I(\alpha,\beta,[C])$.

The key nontrivial property of the ECH index is that if
$C\in\mc{M}(\alpha,\beta)$ is not multiply covered, then
\begin{equation}
\label{eqn:IFact1}
\op{ind}(C) \le I(C) - 2 \delta(C).
\end{equation}
Here $\op{ind}(C)$ denotes the Fredholm index of $C$; if $J$ is
generic, then $\mc{M}(\alpha,\beta)$ is a manifold near $C$ of
dimension $\op{ind}(C)$, see \cite{dragnev}.  Also, $\delta(C)$ is a
nonnegative integer which is zero if and only if $C$ is embedded.  The
index inequality \eqref{eqn:IFact1} was proved in a simpler setting in
\cite[Thm.\ 1.7]{pfh2}, and is proved in the present setting in
\cite[Thm.\ 4.15]{ir} with the help of \cite{siefring1}.

Another useful fact, which is a special case of \cite[Thm.\ 5.1]{ir},
is that if $C$ and $C'$ are two holomorphic curves whose images in
$\R\times Y$ do not have a common irreducible component, then
\begin{equation}
\label{eqn:IFact2}
I(C\cup C') \ge I(C) + I(C') + 2 C\cdot C'.
\end{equation}
Here $C\cdot C'$ denotes the algebraic intersection number of $C$ and
$C'$, a nonnegative integer which is zero if and only if $C$ and $C'$
are disjoint.

The above two inequalities imply the following classification of
(possibly multiply covered) holomorphic curves with low ECH index when
the almost complex structure $J$ is generic.

\begin{notation}
  Any holomorphic curve $C\in\mc{M}(\alpha,\beta)$ can be uniquely
  written as $C=C_0\cup C_1$, where $C_0$ and $C_1$ are holomorphic
  curves such that the image of $C_0$ is a union of $\R$-invariant
  cylinders in $\R\times Y$, while no component of $C_1$ maps to an
  $\R$-invariant cylinder.
\end{notation}

\begin{lemma}
\label{lem:classification}
Suppose $J$ is generic, let $\alpha$ and $\beta$ be ECH generators,
and let $C\in\mc{M}(\alpha,\beta)$.  Write $C=C_0\cup C_1$ as
above. Then:
\begin{description}
\item{(a)} $I(C) \ge 0$, with equality if and only if $C=C_0$.
\item{(b)} If $I(C)\in\{1,2\}$, then $C_0$ and $C_1$ are disjoint in
  $\R\times Y$; $C_1$ is embedded in $\R\times Y$; and
  $\op{ind}(C_1)=I(C_1)=I(C)$.
\end{description}
\end{lemma}

\begin{proof}
  This follows from the inequalities \eqref{eqn:IFact1} and
  \eqref{eqn:IFact2}, as explained eg in \cite[Prop.\ 7.15]{obg1}.
\end{proof}

\subsection{The differential}

The differential
\[
\partial : C_*(Y,\lambda,\Gamma) \longrightarrow
C_{*-1}(Y,\lambda,\Gamma)
\]
is defined as follows.  Fix a generic almost complex structure $J$.
If $\alpha$ is an ECH generator with $[\alpha]=\Gamma$, define
\[
\partial\alpha \eqdef \sum_\beta
\sum_{\substack{C\in\mc{M}(\alpha,\beta)/\R\\I(C)=1}}\epsilon(C) \cdot \beta.
\]
Here the first sum is over ECH generators $\beta$ with
$[\beta]=\Gamma$.  In the second sum, two curves $C=C_0\cup C_1$ and
$C'=C_0'\cup C_1'$ in $\mc{M}(\alpha,\beta)$ are considered equivalent
whenever $C_1=C_1'$.  The $\R$ action on $\mc{M}(\alpha,\beta)$ is
given by translation of the $\R$ coordinate on $\R\times Y$.  Finally,
$\epsilon(C)\in\{\pm 1\}$ is a sign, which depends on some additional
choices described below.  However making different choices to define
the signs will result in isomorphic chain complexes.  It is shown in
\cite[Lem.\ 7.19 and Thm.\ 7.20]{obg1} that $\partial$ is well defined
and $\partial^2=0$.

To determine the signs $\epsilon(C)$, one orients all moduli spaces of
non-multiply covered holomorphic curves (with the ends at positive
hyperbolic orbits ordered and with no end of even multiplicity at a
negative hyperbolic orbit\footnote{For holomorphic curves with even
  multiplicity ends at negative hyperbolic orbits, one needs to
  further choose ``asymptotic markings'' of such ends in order to
  orient the moduli spaces.}) by choosing ``coherent orientations'' of
the relevant deformation operators with the conventions in
\cite[\S9]{obg2}.  The quotients of such moduli spaces by the $\R$
action are then oriented using the ``$\R$-direction first''
convention.  One also fixes an ordering of the set of all positive
hyperbolic embedded Reeb orbits.  Finally, given a curve $C=C_0\cup
C_1$ with $I(C)=1$, writing $C_1\in\mc{M}(\alpha',\beta')$, we define
$\epsilon(C)$ to be the orientation of the point
$C_1\in\mc{M}(\alpha',\beta')/\R$.  This orientation is well defined
thanks to our assumption that no hyperbolic orbit appears in an ECH
generator with multiplicity greater than one.

We denote the homology of this chain complex by
$ECH_*(Y,\lambda,\Gamma)$.  Although the differential $\partial$
depends on $J$, it turns out that the homology of the chain complex
does not.  This follows from the comparison with Seiberg-Witten Floer
homology below.  We also expect that one could prove this directly with
holomorphic curves by a generalization of the proof in
\cite{obg1,obg2} that $\partial^2=0$.

\subsection{The U map}

There is also a degree $-2$ chain map
\[
U:C_*(Y,\lambda,\Gamma) \longrightarrow C_{*-2}(Y,\lambda,\Gamma).
\]
The definition of $U$ was sketched in \cite[\S12]{t3}, and we give
more details here.

To define $U$, fix a generic almost complex structure $J$ and make the
choices described above that are needed to define the signs in the
differential $\partial$.  Also fix a point $y\in Y$ which is not on
any Reeb orbit.  If $\alpha$ and $\beta$ are ECH generators, define
$\mc{M}_y(\alpha,\beta)$ to be the moduli space of pairs $(u,z)$ where
$u:(C,j)\to(\R\times Y,J)$ is an element of $\mc{M}(\alpha,\beta)$,
and $z\in C$ is a marked point with $u(z)=(0,y)$.  Finally, if
$\alpha$ is an ECH generator, define
\[
U\alpha\eqdef
\sum_\beta\sum_{\substack{C\in\mc{M}_y(\alpha,\beta)\\I(C)=2}}\epsilon(C)
\cdot \beta.
\]

Here the sign $\epsilon(C)$ is defined as follows.  Write $C=C_0\cup
C_1$ as usual, and write $C_1\in\mc{M}(\alpha',\beta')$.  Recall from
Lemma~\ref{lem:classification}(b) that $C_1$ is embedded in $\R\times
Y$, and $\mc{M}(\alpha',\beta')$ is a $2$-dimensional manifold near
$C_1$. Let $(v_1,v_2)$ be a positively oriented basis for
$T_{C_1}\mc{M}(\alpha',\beta')$.  The tangent vectors $v_1,v_2$
determine elements $w_1,w_2$ of the normal bundle to $C_1$ at $(0,y)$.
A standard transversality argument shows that if $J$ is generic, then
$w_1,w_2$ are necessarily linearly independent.  We then define
$\epsilon(C)$ to be $+1$ if $(w_1,w_2)$ is a positively oriented basis
for $N_{(0,y)}C_1$, and $-1$ otherwise.

\begin{lemma}
\label{lem:UDefined}
Suppose $J$ is generic.  Then:
\begin{description}
\item{(a)}
$U$ is well defined.
\item{(b)}
Suppose $C\in\mc{M}(\alpha,\beta)$ contributes to $U$, and write
$C=C_0\cup C_1$ as usual.  Then $C_1$ is embedded and connected and
$\op{ind}(C_1)=I(C_1)=2$.
\end{description}
\end{lemma}

\begin{proof}
 Note that if $J$ is generic, then:
\begin{description}
\item{(i)} The elements of $\mc{M}_y(\alpha,\beta)$ with $I=2$ are
  isolated points in $\mc{M}_y(\alpha,\beta)$.
\item{(ii)} $(0,y)$ is not in the image of any holomorphic curve with $I=1$.
\end{description}
Condition (i) follows from the transversality that was discussed in
the previous paragraph.  Condition (ii) for generic $J$ follows from
Lemma~\ref{lem:classification}(b) together with our assumption that
$y$ is not on any Reeb orbit.  Assume that $J$ is generic so that (i)
and (ii) hold.

(a) We need to show that if $\alpha$ and $\beta$ are ECH generators,
then the set of holomorphic curves $C\in\mc{M}_y(\alpha,\beta)$ with
$I(C)=2$ is finite. Suppose to the contrary that $C^{(n)}$ is an
infinite sequence of distinct such curves for $n=1,2,\ldots$.  By
\cite[Lem.\ 7.23]{obg1}, we can replace the sequence $C^{(n)}$ with a
subsequence that converges in the sense of \cite{behwz} to either
a curve in $\mc{M}_y(\alpha,\beta)$, or a ``broken'' curve in which
two of the levels have $I=1$ and the remaining levels have $I=0$.
These cases are impossible by conditions (i) and (ii) above,
respectively.

(b) All except the connectedness follows from
Lemma~\ref{lem:classification}(b).  If $C_1$ is disconnected, then it
follows from Lemma~\ref{lem:classification} and the inequality
\eqref{eqn:IFact2} that $C_1$ has two components and each component
has $I=1$.  This contradicts condition (ii) above.
\end{proof}

One can then show that $U$ is a chain map:
\[
\partial U = U \partial.
\]
The idea of the proof is to count the ends of moduli spaces of $I=3$
curves with marked points mapping to $(0,y)$.  The details are a
straightforward modification of the proof that $\partial^2=0$ in
\cite{obg1,obg2}, because in the analysis, curves with a marked point
constraint behave essentially the same way as curves with index one
less and no marked point constraint.

Moreover, up to chain homotopy, $U$ does not depend on the choice of
$y\in Y$.  To prove this, let $y'\in Y$ be another point not on any
Reeb orbit, and let $U'$ denote the corresponding degree $-2$ chain
map.  Choose a path $\eta$ in $Y$ from $y$ to $y'$ (here we are using
the assumption that $Y$ is connected).  Then counting $I=1$ curves
with marked points mapping to the path $\{0\}\times\eta$ in $\R\times
Y$, using the same signs as in the definition of $\partial$, defines a
degree $-1$ map $K$ such that
\[
\partial K + K \partial = U - U'.
\]
To prove this last equation, one counts ends and boundary points of
the moduli space of $I=2$ curves with marked points mapping to
$\{0\}\times\eta$, using the exact same gluing analysis as in the proof
that $\partial^2=0$.

In conclusion, we obtain a well-defined map on homology
\begin{equation}
\label{eqn:U}
U: ECH_*(Y,\lambda,\Gamma) \longrightarrow ECH_{*-2}(Y,\lambda,\Gamma)
\end{equation}
Again, the comparison with Seiberg-Witten theory shows that this does
not depend on $J$, and we expect that this can also be shown directly
using holomorphic curves.

\subsection{Relation with Seiberg-Witten Floer homology}

It is shown in \cite{echswf} that embedded contact homology is
isomorphic to a version of Seiberg-Witten Floer homology as defined by
Kronheimer-Mrowka \cite{km}.  The precise statement is that for each
$\Gamma\in H_1(Y)$, there is an isomorphism
\begin{equation}
\label{eqn:echswf}
ECH_*(Y,\lambda,\Gamma) \simeq \check{HM}_*(-Y,\frak{s}(\xi) +
\op{PD}(\Gamma))
\end{equation}
of relatively $\Z/d(c_1(\xi)+2\op{PD}(\Gamma))$-graded abelian
groups\footnote{In \cite{echswf}, the right hand side of
  \eqref{eqn:echswf} is replaced by the canonically isomorphic group
  $\widehat{HM}^{-*}(Y,\frak{s}(\xi) + \op{PD}(\Gamma))$.}.  Here
$\frak{s}(\xi)$ denotes the spin-c structure associated to the
oriented $2$-plane field $\xi$ as in \cite[\S28]{km}.

It is further shown in \cite{tu} that, at least up to signs, the
isomorphism \eqref{eqn:echswf} interchanges the map $U$ in
\eqref{eqn:U} with the map
\[
U_\dagger: \check{HM}_*(-Y,\frak{s}(\xi) +
\op{PD}(\Gamma)) \longrightarrow \check{HM}_{*-2}(-Y,\frak{s}(\xi) +
\op{PD}(\Gamma))
\]
defined in \cite{km}.

The above equivalence, together with known properties of $\check{HM}$,
implies the following facts about ECH which we will need.

\begin{proposition}[finiteness]
\label{prop:finite}
\begin{description}
\item{(a)}
For each $\Gamma\in H_1(Y)$ and for each grading $*$, the group
$ECH_*(Y,\lambda,\Gamma)$ is finitely generated.
\item{(b)}
There are only finitely many $\Gamma\in H_1(Y)$ such that
$ECH_*(Y,\lambda,\Gamma)$ is nonzero.
\end{description}
\end{proposition}

\begin{proof}
The corresponding facts about $\check{HM}_*$ are proved in \cite[Lem.\
22.3.3 and Prop.\ 3.1.1]{km}.
\end{proof}

\begin{proposition}[torsion spin-c structures]
\label{prop:torsion}
Let $\Gamma\in H_1(Y)$ and suppose that $c_1(\xi) + 2
\op{PD}(\Gamma)\in H^2(Y;\Z)$ is torsion, so that
$ECH_*(Y,\lambda,\Gamma)$ is relatively $\Z$-graded.  Then:
\begin{description}
\item{(a)}
$ECH_*(Y,\lambda,\Gamma)$ is zero if the grading $*$ is sufficiently small.
\item{(b)}
$ECH_*(Y,\lambda,\Gamma)$ is nonzero for an infinite set of gradings
$*$.
\item{(c)} If the grading $*$ is sufficiently large then the $U$ map
  \eqref{eqn:U} is an isomorphism.
\end{description}
\end{proposition}

\begin{proof}
The corresponding properties of $\check{HM}$ are proved as follows.
Part (a) follows from the definition of $\check{HM}$ in \cite{km}.
Part (b) is proved in \cite[Cor.\ 35.1.4]{km}.  Part (c) is a
consequence of \cite[Lem.\ 33.3.9]{km}.
\end{proof}

\section{More about the holomorphic curves in ECH}
\label{sec:subtleties}

We now give some more detailed information which we will need
concerning the structure of the holomorphic curves that contribute to
the $U$ map in ECH.  (Similar results hold for the curves that
contribute to the ECH differential $\partial$, but will not be needed
here.)  Throughout this section fix $(Y,\lambda)$ as in
\S\ref{sec:ECHReview}, and also fix a generic almost complex structure
$J$ on $\R\times Y$ as needed to define ECH.

\subsection{Possible multiplicities of the ends}
\label{sec:PM}

We begin by recalling some restrictions on the multiplicities of the
ends of holomorphic curves that contribute to $U$.

If $\gamma$ is an embedded elliptic Reeb orbit, and if $\tau$ is a
trivialization of $\xi|_\gamma$, then $\tau$ is homotopic to a
trivialization with respect to which the linearized Reeb flow on the
contact planes along $\gamma$ is rotation by angle $e^{2\pi i\theta}$
for some real number $\theta$, which we call the {\em monodromy
  angle\/} of $\gamma$ with respect to $\tau$.  Our standing
assumption that all Reeb orbits are nondegenerate implies that
$\theta$ is irrational.  Changing the trivialization $\tau$ shifts
$\theta$ by an integer.

If $\theta$ is an irrational number, define $S_\theta$ to be the set
of positive integers $q$ such that $\ceil{q'\theta}/q' >
\ceil{q\theta}/q$ for all $q'\in\{1,\ldots,q-1\}$.  That is, $q\in
S_\theta$ if and only if $\theta$ is better approximated from above by
a rational number with denominator $q$ than by a rational number with
any smaller denominator.  The set $S_\theta$ depends only on the
equivalence class of $\theta$ in $\R/\Z$.

\begin{lemma}
\label{lem:multiplicity}
Let $C$ be a holomorphic curve that contributes to $U$, and write
$C=C_0\cup C_1$ as usual.  Let $\gamma$ be an embedded elliptic Reeb
orbit with monodromy angle $\theta$ with respect to some
trivialization.  Then:
\begin{itemize}
\item
If $C_1$ has a positive end at $\gamma$ of multiplicity $m$, then
$m\in S_{-\theta}$.
\item
If $C_1$ has a negative end at $\gamma$ of multiplicity $m$, then
$m\in S_\theta$.
\end{itemize}
\end{lemma}

\begin{proof}
  By Lemma~\ref{lem:UDefined}(b), the curve $C_1$ is not multiply
  covered and has $\op{ind}(C_1)=I(C_1)$. The conclusions of the lemma
  are then part of the necessary conditions for equality in the index
  inequality \eqref{eqn:IFact1}; see \cite[Thm.\ 1.7]{pfh2} or
  \cite[Thm.\ 4.15]{ir}.
\end{proof}

We now show that the allowable multiplicities have ``density zero''.

\begin{notation}
If $A$ is a subset of the positive integers, define the {\em
  density\/} of $A$ to be
\[
d(A) \eqdef \lim_{N\to\infty}\frac{1}{N}\big|A\cap\{1,\ldots,N\}\big|,
\]
if this limit exists.
\end{notation}

\begin{lemma}
\label{lem:density}
Let $\theta$ be an irrational number.  Then $S_\theta$ has density
zero.
\end{lemma}

\begin{proof}
  Write the elements of $S_\theta$ in increasing order as
  $q_1,q_2,\dots$.  It is enough to show that as $i\to\infty$, the
  differences $q_{i+1}-q_i$ are nondecreasing and tend to infinity.
  In fact these differences are some of the elements of $S_{-\theta}$,
  in increasing order, each repeated some finite number of times.  One
  can prove this by noting that the fractions $\ceil{q\theta}/q$ for
  $q\in S_\theta$ and $\floor{q\theta}/q$ for $q\in S_{-\theta}$ are
  the semiconvergents in the continued fraction expansion of $\theta$,
  and using some basic facts about continued fractions.
\end{proof}

\subsection{Embeddedness in the 3-manifold}
\label{sec:foliation}

The non-$\R$-invariant components of the holomorphic curves counted by
the $U$ map are embedded in $\R\times Y$.  Using arguments going back
to Hofer-Wysocki-Zehnder \cite{hwz2,hwza} and developed further by
Siefring and Wendl, one can show that under certain circumstances the
projections of these curves to $Y$ are also embeddings, and the
corresponding moduli spaces of holomorphic curves locally give a
foliation of $Y$.  In particular, we will need the following
proposition.  Some more general criteria for 3-dimensional
embeddedness and foliations are discussed in \cite{siefring2,wendl3}.

\begin{proposition}
\label{prop:foliation}
Assume that $J$ is generic.  Let $C$ be a holomorphic curve that
contributes to $U$, and decompose $C=C_0\cup C_1$ as usual.  Suppose
that $C_1$ has genus zero, all ends of $C_1$ are at elliptic Reeb
orbits, and $C_1$ does not have two positive ends or two negative ends
at the same Reeb orbit.  Then:
\begin{description}
\item{(a)} The projection of $C_1$ to the three-manifold $Y$ is an
  embedding.  Moreover the projections to $Y$ of the holomorphic
  curves in the same moduli space component as $C_1$ give a foliation
  of some subset of $Y$.
\item{(b)}
$C_1$ does not have both a positive end and a negative end at the same
Reeb orbit.
\end{description}
\end{proposition}

\begin{proof}
The proof has seven steps.

{\em Step 1.\/} Let $C$ be a non-$\R$-invariant connected holomorphic
curve.  We begin by recalling the asymptotic behavior of an end of $C$
at an embedded Reeb orbit $\gamma$; for details see
\cite{hwz1,siefring1,obg2}.

By rescaling we may assume that $\gamma$ is parametrized by
$S^1=\R/\Z$.  The {\em asymptotic operator\/} associated to $\gamma$
is the operator
\[
L_\gamma : C^\infty(\gamma^*\xi) \longrightarrow C^\infty(\gamma^*\xi)
\]
defined by
\[
L_\gamma \eqdef J\nabla_t^R,
\]
where $t$ denotes the $S^1$ coodinate and $\nabla^R$ denotes the
connection on $\gamma^*\xi$ given by the linearized Reeb flow.  More
explicitly, choose a complex linear, symplectic trivialization $\tau$ of
$\gamma^*\xi$; then in this trivialization,
\[
L_\gamma = \sqrt{-1}\frac{d}{dt} + S(t)
\]
where $S(t)$ is a symmetric $2\times 2$ matrix.  If $k$ is a positive
integer, let $\gamma^k$ denote the $k$-fold iterate of $\gamma$, ie the
pullback of $\gamma$ to $\R/k\Z$.  Then the asymptotic
operator associated to $\gamma^k$ is given in the above trivialization by
\begin{equation}
\label{eqn:Lk}
L_{\gamma^k} = \sqrt{-1}\frac{d}{d\tilde{t}} + S(\pi(\tilde{t})),
\end{equation}
where $\tilde{t}$ denotes the $\R/k\Z$ coordinate and
$\pi:\R/k\Z\to\R/\Z$ denotes the projection.

Now identify a tubular neighborhood of $\gamma$ with $S^1\times D$,
where $D$ is a disk in $\C$, such that the derivative of this
identification agrees with $\tau$.  A positive end of $C$ at $\gamma$
of multiplicity $k$ is then described by a map
\[
\begin{split}
[R,\infty) \times \R/k\Z &\longrightarrow \R\times S^1\times D,\\
(s,\tilde{t}) &\longmapsto
(s,\pi(\tilde{t}),\varphi(s,\tilde{t})).
\end{split}
\]
Moreover the function $\varphi$ satisfies
\begin{equation}
\label{eqn:PA}
\varphi(s,\tilde{t}) = e^{-\lambda s}(\nu(\tilde{t}) +\rho(s,\tilde{t})),
\end{equation}
where $\nu:S^1\to\C$ is a nonzero eigenfunction of the asymptotic operator
$L_{\gamma^k}$ with eigenvalue $\lambda>0$, while $\rho(s,\tilde{t})$
and all of its derivatives decay exponentially as $s\to\infty$.  A
negative end of $C$ at $\gamma$ is similarly described by
\begin{equation}
\label{eqn:NA}
\varphi(s,\tilde{t})= e^{-\lambda s}(\nu(\tilde{t}) + \rho(s,\tilde{t}))
\end{equation}
for $s\in(-\infty,R]$, where $\nu$ is a nonzero eigenfunction of
$L_{\gamma^k}$ with eigenvalue $\lambda<0$, and the function
$\rho(s,\tilde{t})$ and all of its derivatives decay exponentially as
$s\to-\infty$.

Note that the eigenfunction $\nu$ in \eqref{eqn:PA} or \eqref{eqn:NA}
can never vanish, by equation \eqref{eqn:Lk} and the uniqueness of
solutions to ODE's.  Furthermore, as shown in \cite{hwz2}, if $\gamma$
is elliptic with monodromy $\theta$ with respect to $\tau$, then for a
positive end $\nu$ has winding number at most
\begin{equation}
\label{eqn:PWB}
\op{wind}(\nu) \le \floor{k\theta},
\end{equation}
while for a negative end $\nu$ has winding number at least
\begin{equation}
\label{eqn:NWB}
\op{wind}(\nu) \ge \ceil{k\theta}.
\end{equation}

{\em Step 2.\/} We now deduce an important inequality.  Suppose that
the holomorphic curve $C$ is connected, immersed, and
non-$\R$-invariant.  Let $\psi$ denote the section of the normal
bundle to $C$ given by the projection of $\partial/\partial s$, where
$s$ denotes the $\R$ coordinate on $\R\times Y$.  A calculation using
the above asymptotic formulas and winding bounds along with similar
winding bounds for the hyperbolic ends, cf.\ \cite{hwz2}, shows that
the algebraic count of zeroes of $\psi$ is finite and satisfies
\begin{equation}
\label{eqn:noZeroes}
2\#\psi^{-1}(0) \le 2g(C) - 2 + \op{ind}(C) + h_+(C),
\end{equation}
where $h_+(C)$ denotes the number of ends of $C$ at positive
hyperbolic orbits (or at negative hyperbolic orbits with even
multiplicity).

{\em Step 3.\/} Now write $C\eqdef C_1$. Lemma~\ref{lem:UDefined}
implies that $C$ is connected, embedded in $\R\times Y$, and has
Fredholm index $\op{ind}(C) = 2$.  Thanks to our hypotheses and the
fact that $\op{ind}(C)=2$, the right hand side of \eqref{eqn:noZeroes}
equals zero.  On the other hand, because $\psi$ satisfies a linear PDE
with the same symbol as a Cauchy-Riemann equation, all zeroes of
$\psi$ have positive multiplicity.  Consequently $\psi$ has no zeroes,
and it follows that the projection of $C$ to $Y$ is at least an
immersion.

{\em Step 4.\/} Let $\nu$ denote the asymptotic eigenfunction
associated to an end of $C$ at $\gamma$ with
multiplicity $k$.  We claim that there does not exist a solution to
the equation
\[
\nu(\tilde{t}) = \rho \nu(\tilde{t} + l), \quad\quad
\rho>0,\;\;\tilde{t}\in\R/k\Z, \;\; l\in\{1,\ldots,k-1\}.
\]

If a solution to the above equation exists, then it follows from
\eqref{eqn:Lk} and the uniqueness of solutions to ODE's that $\rho=1$
and $\nu$ is the pullback of an eigenfunction of $L_{\gamma^{k'}}$
where $k'<k$ is a divisor of $k$.  But this cannot happen because the
winding number $\op{wind}(\nu)$ is relatively prime to $k$.  To prove
this last fact, one notes that the winding bounds \eqref{eqn:PWB} and
\eqref{eqn:NWB} are sharp here, eg because equality holds in
\eqref{eqn:noZeroes}, and then applies Lemma~\ref{lem:multiplicity}.

{\em Step 5.\/} For $t > 0$, let $C(t)$ denote the holomorphic curve
in $\R\times Y$ obtained from $C$ by translation in the $\R$ direction
by distance $t$.

Clam: There exists $R>0$ such that if $t>0$ is
sufficiently small, then any intersection of $C$ with $C(t)$ has
$|s|<R$.

Proof of Claim: Two different ends of $C$ and $C(t)$ cannot intersect
by the last of our hypotheses.  To show that an end of $C$ cannot
intersect the same end of $C(t)$, one uses the asymptotic formulas
\eqref{eqn:PA} and \eqref{eqn:NA} together with Step 4.

It follows from the above claim that $C$ is disjoint from $C(t)$ when
$t>0$ is sufficiently small.  Otherwise we could take a sequence of
intersections $\{x_n\}$ of $C$ with $C(t_n)$ where $t_n\to 0$, and by
the claim we could pass to a convergent subsequence.  The limit of
this subsequence would then be a zero of $\psi$, contradicting Step 3.

{\em Step 6.\/} We now complete the proof of part (a).  To prove that
the projection of $C$ to $Y$ is an embedding, it is enough to show
that $C$ is disjoint from $C(t)$ for all $t>0$.  More generally, to
prove all of part (a) it is enough to show that if $C'\neq C$ is any
holomorphic curve in the same moduli space component as $C$, then $C$
and $C'$ are disjoint.

It follows from the detailed asymptotics in \cite{siefring1} that $C$
and $C'$ have only finitely many intersections.  By \cite[Lem.\
8.5]{pfh2}, the algebraic count of intersections of $C$ and $C'$ is
given by
\[
C\cdot C' = Q_\tau(C) + \ell_\tau(C,C').
\]
Here $\tau$ is a trivialization of $\xi$ over all the Reeb orbits at
which $C$ has ends, and $\ell_\tau(C,C')$ denotes the ``asymptotic
linking number'' of $C$ and $C'$ with respect to $\tau$, defined in
\cite[\S8.2]{pfh2} or \cite[\S2.7]{ir}.

The asymptotic linking number is bounded from above by
\begin{equation}
\label{eqn:linkingBound}
\ell_\tau(C,C') \le \sum_{\gamma\in P_+}k\floor{k\theta} -
\sum_{\gamma \in P_-}k\ceil{k\theta}.
\end{equation}
Here the first sum is over the Reeb orbits $\gamma$ at which $C$ has a
positive end, the second sum is over the Reeb orbits $\gamma$ at which
$C$ has a negative end, and in each summand, $\theta$ denotes the
monodromy angle of $\gamma$ with respect to $\tau$ and $k$ denotes the
multiplicity of the corresponding end of $C$.  The inequality
\eqref{eqn:linkingBound} is a special case of a linking bound which is
proved in a simpler situation in \cite[Lem. 6.9]{pfh2} and which
follows in the present case by the asymptotic analysis in
\cite{siefring1}.  Thus we obtain an upper bound on the algebraic
intersection number,
\begin{equation}
\label{eqn:intersectionBound}
C\cdot C' \le Q_\tau(C) + \sum_{\gamma\in P_+}k\floor{k\theta} -
\sum_{\gamma \in P_-}k\ceil{k\theta}.
\end{equation}

The right hand side of \eqref{eqn:intersectionBound} is a topological
invariant of $C$ which does not depend on $C'$. If $C'=C(t)$ where
$t>0$ is small, then it follows from the asymptotics discussed above
that the inequality \eqref{eqn:linkingBound}, and hence the inequality
\eqref{eqn:intersectionBound}, is sharp.  On the other hand we know
from Step 5 that $C\cdot C'=0$ in this case.  Thus
\eqref{eqn:intersectionBound} says that $C\cdot C'\le 0$ for all $C'$.
It follows by intersection positivity that $C$ and $C'$ are disjoint
for all $C'$.

{\em Step 7.\/} We now prove part (b).  Suppose $\gamma$ is an
embedded elliptic Reeb orbit at which $C$ has both a positive and a
negative end.  Let $T$ denote the boundary of a small tubular
neighborhood of $\gamma$.  Let $\tau$ be a trivialization of $\xi$
over $\gamma$ and identify $T\simeq S^1\times S^1$ compatibly with
this trivialization, where the first $S^1$ factor is identified with
$\gamma$.  Let $\theta$ denote the monodromy angle of $\gamma$ with
respect to $\tau$.

It follows from the above asymptotic formulas and winding bounds that
the projection of the positive end of $C$ to $Y$ intersects $T$
transversely in a circle representing a homology class $(q_+,p_+)\in
H_1(T^2)$, such that $q_+$ is the multiplicity of the end and
$p_+/q_+<\theta$.  Likewise, the projection of the negative end of $C$
to $Y$ intersects $T$ transversely in a circle with homology class
$(q_-,p_-)$ satisfying $p_-/q_->\theta$.  Since $p_+/q_+\neq p_-/q_-$,
the positive and negative circles in $T$ must intersect, contradicting
part (a).
\end{proof}

\subsection{Euler characteristic}

We next recall from \cite[\S6]{ir} a variant of the ECH index, denoted
by $J_0$, which bounds the negative Euler characteristic of holomorphic
curves, similarly to the way that the ECH index $I$ bounds the Fredholm
index in \eqref{eqn:IFact1}.

If $\alpha=\{(\alpha_i,m_i)\}$ and $\beta=\{(\beta_j,n_j)\}$ are ECH
generators with $[\alpha]=[\beta]$, and if $Z\in H_2(Y,\alpha,\beta)$,
one defines
\begin{equation}
\label{eqn:J0Def}
J_0(\alpha,\beta,Z) \eqdef -c_\tau(Z) + Q_\tau(Z) +
\sum_i\sum_{k=1}^{m_i-1}\op{CZ}_\tau(\alpha_i^k) -
\sum_j\sum_{k=1}^{n_j-1}\op{CZ}_\tau(\beta_j^k).
\end{equation}
Here $\tau$ is a trivialization of $\xi$ over the Reeb orbits
$\alpha_i$ and $\beta_j$;  one can check that $J_0$, like $I$, does
not depend on $\tau$, even though the individual terms in its
definition do.  Also, like the ECH index, $J_0$ is additive in the
sense that
\begin{equation}
\label{eqn:J0Additive}
J_0(\alpha,\beta,Z)+J_0(\beta,\gamma,W) = J_0(\alpha,\gamma,Z+W).
\end{equation}
If $C\in\mc{M}(\alpha,\beta)$ is a holomorphic
curve, we write $J_0(C)\eqdef J_0(\alpha,\beta,[C])$.

We now have the following bound on topological complexity in terms of $J_0$.

\begin{lemma}
\label{lem:complexity}
Let $C\in\mc{M}(\alpha,\beta)$ be a holomorphic curve
that contributes to $U$.  Write $C=C_0\cup C_1$ as usual.  
Then
\begin{equation}
\label{eqn:control}
J_0(C) \ge 2g(C_1) - 2 + \sum_i(2n_i^+ + t_i^+ - 1) +
\sum_j(2n_j^- + t_j^- - 1).
\end{equation}
Here $g(C_1)$ denotes the genus of $C_1$; $n_i^+$ denotes the number
of positive ends of $C_1$ at $\alpha_i$, and $n_j^-$ denotes the
number of negative ends of $C_1$ at $\beta_j^-$; $t_i^+$ is $1$ if the
image of $C_0$ contains $\R\times\alpha_i$ and $0$ otherwise; and
$t_j^-$ is $1$ if the image of $C_0$ contains $\R\times\beta_j$ and
$0$ otherwise\footnote{In fact equality holds in \eqref{eqn:control}
  here.  One can show this by the
  arguments in \cite{ir}, or by a more direct calculation using the
  necessary conditions for equality in \eqref{eqn:IFact1}.  However we
will not need this.}.
\end{lemma}

\begin{proof}
  Since $C_1$ is embedded by Lemma~\ref{lem:UDefined}(b), one can apply
  \cite[Prop.\ 6.9]{ir} to obtain
\begin{equation}
\label{eqn:J01}
J_0(C_1) \ge 2g(C_1) - 2 + \sum_{i: n_i^+>0}(2n_i^+ - 1) +
\sum_{j:n_j^->0}(2n_j^- - 1).
\end{equation}
Also
\cite[Prop.\ 6.14]{ir} implies that
\begin{equation}
\label{eqn:J02}
J_0(C) \ge J_0(C_0) + J_0(C_1) + \sum_{i: n_i^+>0}t_i^+ +
\sum_{j:n_j^->0}t_j^-.
\end{equation}
Finally, it follows from the definition of $J_0$ that
\begin{equation}
\label{eqn:J03}
J_0(C_0)=0.
\end{equation}
Combining \eqref{eqn:J01}, \eqref{eqn:J02}, and \eqref{eqn:J03} proves
the lemma.
\end{proof}

It follows immediately from the inequality \eqref{eqn:control} that
\begin{equation}
\label{eqn:J0ge-1}
J_0(C) \ge -1,
\end{equation}
since all of the summands in the sums over $i$ and $j$ are nonnegative
and $C_1$ has at least one end.  We also deduce the following
criterion for recognizing holomorphic cylinders in a certain situation
which will arise later.

\begin{lemma}
\label{lem:control}
Let $\gamma_1$ and $\gamma_2$ be distinct embedded elliptic Reeb
orbits and let $\alpha=\gamma_1^{m_1}\gamma_2^{m_2}$ and
$\alpha'=\gamma_1^{m_1'}\gamma_2^{m_2'}$ with $m_1,m_2,m_1',m_2'\neq 0$.
Suppose $C\in\mc{M}(\alpha,\alpha')$ contributes to $U$, and write
$C=C_0\cup C_1$ as usual.  Assume that $C_1$ has ends at both
$\gamma_1$ and $\gamma_2$. Then:
\begin{description}
\item{(a)} $J_0(C) \ge 2$.
\item{(b)} If $J_0(C)=2$, then $C_1$ is a cylinder.
\end{description}
\end{lemma}

\begin{proof}
 (a) In this situation the inequality \eqref{eqn:control} can be
  rewritten as
\begin{equation}
\label{eqn:control2}
\frac{1}{2}J_0(C) \ge g(C_1) - 3 + \sum_{i=1}^2\left(n_i^+ + n_i^- +
    T_i\right),
\end{equation}
where $T_i$ is defined to be $1$ if the image of $C_0$ contains
$\R\times\gamma_i$ and $0$ otherwise.  Since $C_1$ has ends at both
$\gamma_1$ and $\gamma_2$, we have $n_i^+ + n_i^-\ge 1$.  Also, since
$m_i,m_i'\neq 0$, it follows that
\begin{equation}
\label{eqn:control3}
T_i=0 \Longrightarrow n_i^+,n_i^-\ge 1.
\end{equation}
So in all cases we have
\begin{equation}
\label{eqn:control4}
n_i^+ + n_i^- + T_i \ge 2.
\end{equation}
Putting \eqref{eqn:control4} into \eqref{eqn:control2} gives
$J_0(C)\ge 2$.

(b) By the above, if $J_0(C)=2$ then $g(C_1)=0$ and $n_i^+ + n_i^- +
T_i= 2$ for $i=1,2$.  If $T_1=T_2=1$, then it follows immediately that
$C_1$ is a cylinder so we are done.  If some $T_i=0$, then
$n_i^+=n_i^-=1$ by \eqref{eqn:control3}.  But this contradicts
Proposition~\ref{prop:foliation}(b).
\end{proof}

\section{Contact 3-manifolds with all Reeb orbits elliptic}
\label{sec:lens}

We now prove Theorem~\ref{thm:lens}.  The proof occupies
\S\ref{sec:BO}-\ref{sec:RSA} below.  We then prove
Theorem~\ref{thm:three} in \S\ref{sec:three}.

\subsection{Initial input from Seiberg-Witten theory}
\label{sec:BO}

Throughout the proof of Theorem~\ref{thm:lens}, 
fix a closed oriented connected 3-manifold $Y$ with a contact form
$\lambda$ such that all Reeb orbits are nondegenerate and elliptic.
Also fix a generic almost complex structure on $\R\times Y$ in order to
define the ECH chain complex.

We begin by using Propositions \ref{prop:finite} and
\ref{prop:torsion} from Seiberg-Witten theory to make some basic
observations which will be used repeatedly below.

\begin{lemma}
\label{lem:basic}
Under our assumption that all Reeb orbits are nondegenerate and elliptic:
\begin{description}
\item{(a)}
All differentials in the ECH chain complex vanish.
\item{(b)}
If $\gamma$ is a Reeb orbit then $[\gamma]$ is torsion in $H_1(Y)$.
\item{(c)}
$c_1(\xi)$ is torsion in $H^2(Y;\Z)$.
\end{description}
\end{lemma}

\begin{proof}
  (a) Equation \eqref{eqn:I2} and our assumption that all Reeb orbits
  are elliptic imply that the relative index $I(\alpha,\beta)$ is
  always even, so the differential vanishes.

  (b) If $\gamma$ is an embedded Reeb orbit, then for each nonnegative
  integer $m$, the ECH generator $\gamma^m$ represents a nonzero class
  in $ECH_*(Y,\lambda,m[\gamma])$.  If $[\gamma]$ is not torsion, then
  the homology classes $\{m[\gamma]\}_{m=0,1,\ldots}$ give infinitely
  many $\Gamma$ for which $ECH_*(Y,\lambda,\Gamma)$ is nonzero,
  contradicting Proposition~\ref{prop:finite}(b).

  (c) Since $c_1(\xi)$ is divisible by $2$ in $H^2(Y;\Z)$, there
  exists $\Gamma\in H_1(Y)$ such that $c_1(\xi)+2\op{PD}(\Gamma)=0$.
  By Proposition~\ref{prop:torsion}(b), there exists an admissible
  orbit set $\alpha=\{(\alpha_i,m_i)\}$ with
  $[\alpha]=\sum_im_i[\alpha_i]=\Gamma$.  It then follows from part
  (b) that $c_1(\xi)$ is torsion.
\end{proof}

This lemma simplifies the computation of the ECH index as follows.
Let $\alpha$ and $\beta$ be two ECH generators with
$[\alpha]=[\beta]\in H_1(Y)$, let $\tau$ be a trivialization of $\xi$
over the Reeb orbits in $\alpha$ and $\beta$, and let $Z\in
H_2(Y,\alpha,\beta)$.  Then by equation \eqref{eqn:cZZ'}, the relative
first Chern class $c_\tau(Z)$ depends only on $\alpha$, $\beta$, and
$\tau$, so we can denote it by $c_\tau(\alpha,\beta)$.  Likewise, by
equation \eqref{eqn:QZZ'}, the relative intersection pairing
$Q_\tau(Z)$ depends only on $\alpha$, $\beta$, and $\tau$, so we can
denote it by $Q_\tau(\alpha,\beta)$.

\subsection{Computing the ECH index}
\label{sec:computeI}

We now compute the ECH index of all of the ECH generators.  For
simplicity, let us temporarily assume that all Reeb orbits are
nullhomologous; we will remove this assumption in \S\ref{sec:RSA}.

Let $\gamma_1,\ldots,\gamma_n$ denote the distinct embedded Reeb
orbits.  (At this point in the argument there could be infinitely many
of them, in which case $n$ should be replaced by $\infty$ in the
summations below.)  The ECH generators have the form
$\gamma_1^{m_1}\cdots\gamma_n^{m_n}$ where $m_1,\ldots,m_n$ are
nonnegative integers.  This product notation is shorthand for the
orbit set $\{(\gamma_i,m_i)\mid i=1,\ldots,n; \; m_i\neq 0\}$.

Since $c_1(\xi)$ is torsion and all of the Reeb orbits are
nullhomologous, the relative index on ECH has a unique refinement to
an absolute index which assigns to each generator $\alpha$ an integer
$I(\alpha)$ such that $I(\alpha,\beta)=I(\alpha)-I(\beta)$ and
$I(\emptyset)=0$.  To describe this integer, fix a trivialization
$\tau$ of the contact structure $\xi$ over the $\gamma_i$'s.  The
index of an ECH generator $\alpha=\gamma_1^{m_1}\cdots\gamma_n^{m_n}$
is then given by
\[
I(\alpha) = c_\tau(\alpha) + Q_\tau(\alpha) + \sum_{i=1}^n
\sum_{k=1}^{m_i} \op{CZ}_\tau(\gamma_i^k).
\]
Here $c_\tau(\alpha)$ is shorthand for $c_\tau(\alpha,\emptyset)$, and
$Q_\tau(\alpha)$ is shorthand for $Q_\tau(\alpha,\emptyset)$.

To make this more explicit, define $c_i\eqdef
c_\tau(\gamma_i,\emptyset)$ and $Q_i\eqdef
Q_\tau(\gamma_i,\emptyset)$.  Also, for $i\neq j$ define $Q_{ij}\eqdef
Q_\tau(Z_i,Z_j)$ where $Z_i\in H_2(Y,\gamma_i,\emptyset)$ and $Z_j\in
H_2(Y,\gamma_j,\emptyset)$.  In fact it follows from the definition of
$Q$ in \cite{pfh2,ir} that $Q_{ij}$ does not depend on $Z_i$, $Z_j$,
or $\tau$ and is just the linking number of $\gamma_i$ and $\gamma_j$.
Finally, let $\theta_i$ denote the monodromy angle of $\gamma_i$ with
respect to $\tau$, see \S\ref{sec:PM}.  Since $c_\tau$ is linear in
the relative homology class, we have
\[
c_\tau(\alpha) = \sum_{i=1}^n m_i c_i.
\]
Also $Q_\tau$ is quadratic in the sense that
\begin{equation}
\label{eqn:quadratic}
Q_\tau(\alpha) = \sum_{i=1}^n m_i^2 Q_i + \sum_{i\neq j}m_im_jQ_{ij},
\end{equation}
see \cite[Eq.\ (68)]{pfh2} or \cite[Eq.\ (3.11)]{ir}.
Finally, the Conley-Zehnder terms are given explicitly by
\[
\op{CZ}_\tau(\gamma_i^k) =
2\floor{k\theta_i} + 1.
\]
To simplify the resulting expression for the ECH index, define
\[
\eta_i \eqdef (c_i - Q_i + 1)/2, \quad \quad \phi_i \eqdef Q_i +
\theta_i.
\]
We then obtain
\begin{equation}
\label{eqn:I}
\frac{1}{2}I(\alpha) = \sum_{i=1}^n m_i\eta_i + \sum_{i < j}m_im_j Q_{ij} +
\sum_{i=1}^n \sum_{k=1}^{m_i}\floor{k\phi_i}. 
\end{equation}
We remark that the quantities $\eta_i$ and $\phi_i$ are natural to
consider because they do not depend on $\tau$ (even though $c_i$,
$Q_i$, and $\theta_i$ do).  In fact $Q_i-c_i$ agrees with a familiar
quantity from contact topology, namely the self-linking number of the
transverse knot $\gamma_i$, see eg \cite[\S3.5.2]{geiges}; while
$\phi_i$ is some irrational number.

\begin{example}
\label{ex:ellipsoid}
Suppose $Y$ is an ellipsoid
\[
\left\{(z_1,z_2) \,\big|\, \frac{|z_1^2|}{a_1} +
  \frac{|z_2|^2}{a_2}=1\right\} \subset \C^2,
\]
where $a_1/a_2$ is irrational, with the standard contact form.  Here
there are exactly two Reeb orbits, both elliptic, given by the circles
$z_1=0$ and $z_2=0$.  One can calculate that
$\eta_1=\eta_2=Q_{12}=\phi_1\phi_2=1$.  It is an exercise to deduce
from equation \eqref{eqn:I} that there is exactly one generator of
each nonnegative even index\footnote{The solution to the exercise is
  to consider the line $L$ in the plane with slope $-\phi_1$ passing
  through the point $(m_1,m_2)$.  The right side of \eqref{eqn:I} is
  then the number of lattice points in the triangle consisting of the
  line $L$ and the coordinate axes (including lattice points on the
  boundary), minus $1$.  As one moves the line $L$ up and to the
  right, keeping its slope fixed, one hits all of the lattice points
  in the positive quadrant in succession, and each lattice point has
  index $2$ greater than the previous one.}, as there should be since
$Y\simeq S^3$ and the Seiberg-Witten Floer homology $\check{HM}$ of
$S^3$ is known to have one generator in each nonnegative even degree.
It is interesting to compare this example with the characterization of
ellipsoids in terms of linearized contact homology in \cite{bce}.
\end{example}

\subsection{Proof that there are exactly two embedded orbits}
\label{sec:n=2}

We now prove that $n=2$.  The idea is to use
Proposition~\ref{prop:torsion} and argue that if $n=1$ then there are
not enough generators in a given index range, and if $n\ge 3$ then
there are too many.

By estimating $k\phi_i - 1 \le \floor{k\phi_i}\le k\phi_i$ in
equation \eqref{eqn:I}, we find that the ECH index is approximated by
\begin{equation}
\label{eqn:IApprox}
I(\alpha) = \overline{Q}(m_1,\ldots,m_n)  +
O(m_1+\cdots+m_n),
\end{equation}
where $\overline{Q}$ denotes the quadratic form
\[
\overline{Q}(m_1,\ldots,m_n) \eqdef
 \sum_{i=1}^n m_i^2\phi_i + \sum_{i\neq j}m_im_j Q_{ij}.
\]
Also recall that the numbers $\phi_i$ are irrational.  It follows that
$n\ge 2$, because if $n<2$, then the number of ECH generators with
index $\le k$ is bounded from above by a linear function of $\sqrt{k}$
plus a constant, so there are not enough generators of large index to
satisfy Proposition~\ref{prop:torsion}(b),(c).

On the other hand, it follows from \eqref{eqn:IApprox} that there is a
constant $c$ such that
\[
I(\alpha) \le c(m_1^2+\cdots+m_n^2+1).
\]
Hence the number of ECH generators with index $\le k$ is bounded from
below by a constant times $(k-c)^{n/2}$.  But by
Propositions~\ref{prop:finite}(a) and~\ref{prop:torsion}(a),(c), the
number of generators of index $\le k$ is bounded from above by a
linear function of $k$ plus a constant.  Thus $n\le 2$.

\subsection{A lower bound on the ECH index}
\label{sec:LBI}

With $n=2$ proved, we now establish an additional estimate on the ECH
index which will be needed later:

\begin{lemma}
\label{lem:LBI}
There are constants $c_1,c_2>0$ such that if
$\alpha=\gamma_1^{m_1}\gamma_2^{m_2}$ with $m_1,m_2\ge 0$, then
\[
I(\alpha) \ge c_1(m_1^2+m_2^2)-c_2.
\]
\end{lemma}

\begin{proof}
  By the estimate \eqref{eqn:IApprox}, if
  $\alpha=\gamma_1^{m_1}\gamma_2^{m_2}$ with $m_1,m_2\ge 0$ then
\begin{align}
\label{eqn:QL1}
I(\alpha) &\ge \overline{Q}(m_1,m_2) - L(m_1,m_2),\\
\label{eqn:QL2}
I(\alpha) &\le
\overline{Q}(m_1,m_2)+L(m_1,m_2),
\end{align}
where $L$ is a linear function.  By the lower bound in
\eqref{eqn:QL1}, it is enough to show that $\overline{Q}(m_1,m_2)>0$
whenever $m_1$ and $m_2$ are nonnegative and not both zero.  It
follows from Proposition~\ref{prop:torsion}(a) and the upper bound in
\eqref{eqn:QL2} that $\overline{Q}$ is nonnegative on all lines of
rational slope in the quadrant $\{(m_1,m_2)\mid m_1,m_2\ge 0\}$.
Therefore $\overline{Q}$ is nonnegative on the whole quadrant.
Moreover $\overline{Q}$ is positive on the coordinate axes (minus the
origin), because $\phi_1$ and $\phi_2$ are irrational.  So we just
need to rule out the case where $\overline{Q}$ is degenerate and its
null space has positive slope.

In the case to be ruled out,
\[
\overline{Q}(m_1,m_2) = \phi_1(m_1-\lambda m_2)^2
\]
where $\phi_1,\lambda>0$.  We can then complete the square in the upper bound
\eqref{eqn:QL2} to obtain
\begin{equation}
\label{eqn:CS}
I(\alpha) \le \phi_1(m_1-\lambda m_2 + a_1)^2 + a_2m_1 + a_3
\end{equation}
where $a_1,a_2,a_3$ are constants.  But this contradicts the linear
growth in the number of ECH generators with index $\le k$ as
$k\to\infty$.  To see this, note that if $a_2\le 0$, then by taking
lattice points near the line $m_1=\lambda m_2$ we can find infinitely
many ECH generators with index bounded from above.  And if $a_2>0$,
then calculation using \eqref{eqn:CS} shows that the number of ECH
generators with index $\le k$ grows as at least $k^{3/2}$.
\end{proof}

We remark that in Example~\ref{ex:ellipsoid}, the quadratic form
$\overline{Q}$ is degenerate, and its null space has negative slope.

\subsection{Existence of a holomorphic cylinder}

We now put everything together to prove:

\begin{lemma}
\label{lem:prob}
Let $\gamma_1$ and $\gamma_2$ denote the distinct embedded Reeb orbits
as above.  Then there exist positive integers $m_1,m_2$ and an $I=2$
cylinder $C$ in $\mc{M}(\gamma_1^{m_1},\gamma_2^{m_2})$ or
$\mc{M}(\gamma_2^{m_2},\gamma_1^{m_1})$ or
$\mc{M}(\gamma_1^{m_1}\gamma_2^{m_2},\emptyset)$.
\end{lemma}

\begin{proof}
  By Proposition~\ref{prop:torsion}, if $k_0$ is a sufficiently large
  integer, then there exist ECH generators of index $2k_0$, and the
  $U$ map is an isomorphism in all degrees higher than $2k_0$.  Let
  $\alpha(0)$ be a generator with $I(\alpha(0))=2k_0$.  By induction
  on $k$ we can find a sequence of ECH generators
\[
\left\{\alpha(k) = \gamma_1^{m_1(k)}\gamma_2^{m_2(k)}\right\}
\]
indexed by $k\ge 0$ such that for
  each $k\ge 1$, there exists
\[
C(k)\in\mc{M}(\alpha(k),\alpha(k-1))
\]
which contributes to $U$.  In particular,
\begin{equation}
\label{eqn:Iak}
I(\alpha(k)) = 2(k_0+k).
\end{equation}

Now we consider positive integers $k$ such that the curves $C(k)$ have
certain desirable properties.  Let
\begin{align*}
  A_1 & \eqdef \{\mbox{$k>0\mid m_i(k), m_i(k-1)\neq 0$ for $i=1,2$}\},\\
  A_2 & \eqdef \{\mbox{$k>0\mid m_i(k)\notin S_{-\theta_i}$ and $m_i(k-1)\notin
  S_{\theta_i}$ for $i=1,2$}\},\\
  A_3 & \eqdef \{k>0\mid \mbox{$C(k)_1$ has ends at both $\gamma_1$
    and $\gamma_2$, or $J_0(C(k))\ge 3$}\},\\
  A_4 & \eqdef \{k>0\mid J_0(C(k)) = 2\}.
\end{align*}
Here in the definition of $A_3$ we are decomposing $C(k)=C(k)_0\cup
C(k)_1$ as usual.  It is enough to show that $A_1\cap A_3\cap A_4$ is
nonempty, because if $k\in A_1\cap A_3\cap A_4$ then
Lemma~\ref{lem:control}(b) is applicable to the curve $C(k)$, and
$C(k)_1$ is the cylinder we are seeking.  To prove that $A_1\cap
A_3\cap A_4$ is nonempty, we will show that for each $j=1,2,3,4$ the
set $A_j$ has density $1$.  (The fact that $A_2$ has density $1$ will
be used in the proofs that $A_3$ and $A_4$ do.)

To show that $A_1$ has density $1$, note that since $I(\alpha(k))$
grows linearly with $k$, see equation \eqref{eqn:Iak}, it follows from
Lemma~\ref{lem:LBI} that
\begin{description}
\item{(*)}
There is a constant $c$ such that for each positive integer $N$,
the $N+1$ points
$
(m_1(0),m_2(0)), \ldots,
(m_1(N),m_2(N))
$
are all contained in a ball of radius $c(\sqrt{N}+1)$ centered at the origin.
\end{description}
Consequently $|\{1,\ldots,N\}\setminus A_1|$ grows as at most
$\sqrt{N}$, so $A_1$ has density $1$.  Likewise, (*) and
Lemma~\ref{lem:density} imply that $A_2$ has density $1$.

We now show that $A_3$ has density $1$.  Since $A_1$ and $A_2$ have
density $1$, it is enough to show that $(A_1\cap A_2) \setminus A_3$
has density $0$.  Suppose $k\in (A_1\cap A_2)\setminus A_3$ and let
$C\eqdef C(k)$.  Since $k\notin A_3$, without loss of generality $C_1$
has ends only at $\gamma_1$, and $J_0(C) \le 2$.  These two
conditions, together with the assumption $k\in A_1$ and
Lemma~\ref{lem:complexity}, imply that
\[
g(C_1) + n_1^+ + n_1^- + T_1 \le 3.
\]
Here $n_1^+$ denotes the number of positive ends of $C_1$ at
$\gamma_1$; $n_1^-$ denotes the number of negative ends of $C_1$ at
$\gamma_1$; and $T_1$ is defined to be $1$ if the image of $C_0$
contains $\R\times\gamma_1$ and zero otherwise. Note also that
$n_1^+\ge 1$ since any holomorphic curve of the type we are
considering must have at least one positive end.  We then deduce from
the above inequality that at least one of the following cases holds:
\begin{description}
\item{(i)} $C_1$ has no negative end.
\item{(ii)} 
The image of $C_0$ does not contain
  $\R\times\gamma_1$, and
$C_1$ has exactly one positive end or exactly one negative end.
\item{(iii)} $C_1$ is a cylinder with both positive and negative ends.
\end{description}
Now case (ii) is impossible by our assumption that $k\in A_2$ and
Lemma~\ref{lem:multiplicity}.  And case (iii) is impossible by
Proposition~\ref{prop:foliation}(b).  So to complete the proof that
$A_3$ has density $1$, it is enough to show that case (i) only happens
for $k$ in a set of density zero.

If case (i) holds, let $\widetilde{C}_1$ denote the union of $C_1$
with the part of $C_0$ that maps to $\R\times\gamma_1$.  By the
superadditivity of the ECH index in equation \eqref{eqn:IFact2} and
Lemma~\ref{lem:UDefined}(b), we have
$
2=I(C) \ge I(\widetilde{C}_1) \ge I(C_1)=2,
$
so
\[
I(\widetilde{C}_1)=I(C_1)=2.
\]
That is, writing $m\eqdef m_1(k)$
and $m'\eqdef m_1(k-1)$, we have
\[
I(\gamma_1^{m}) - I(\gamma_1^{m'}) = I(\gamma_1^{m-m'}) = 2.
\]
By equation \eqref{eqn:I}, this is equivalent to
\[
\sum_{k=1}^{m-m'} \floor{(m'+k)\phi_1} = \sum_{k=1}^{m-m'}\floor{k\phi_1} =
1 - \eta_1(m-m').
\]
Recall that $\phi_1>0$.  Now the left equality requires that $m'\le
\ceil{1/\phi_1}$, because otherwise each term in the left sum would be
greater than the corresponding term in the right sum.  And the right
equality can only hold for finitely many values of $m-m'$, because the
left side is approximated by a quadratic function of $m-m'$ while the
right side is a linear function of $m-m'$.  We conclude that case (i)
can only hold for finitely many pairs $(m_1(k),m_1(k-1))$, and by (*)
again this can only happen for $k$ in a set of density zero.

To show that $A_4$ has density $1$, we first show that $J_0$ is
``close'' to $I$.  In the present situation, the relative index $J_0$,
just like $I$, can be uniquely refined to an absolute index which
associates to each ECH generator an integer, such that
$J_0(\alpha,\beta)=J_0(\alpha)-J_0(\beta)$ and $J_0(\emptyset)=0$.
Similarly to \eqref{eqn:I}, we find that if
$\alpha=\gamma_1^{m_1}\gamma_2^{m_2}$, then
\[
\frac{1}{2}J_0(\alpha) = \sum_{i=1}^2m_i(1-\eta_i) + m_1m_2Q_{12} +
\sum_{i=1}^2\sum_{k=1}^{m_i-1}\floor{k\phi_i} - \frac{1}{2}\#\{i\mid
m_i\neq 0\}.
\]
Subtracting this from equation \eqref{eqn:I}, we obtain
\[
I(\alpha) - J_0(\alpha) = \sum_{i=1}^2 m_i(4\eta_i-2) + 2\sum_{i=1}^2
\floor{m_i\phi_i} + \#\{i\mid m_i\neq 0\}.
\]
This equation implies that there is a constant $c$ such that
\[
|J_0(\alpha) - I(\alpha)| \le c(m_1+m_2).
\]
It follows using (*) that for any $\epsilon>0$, if $N$ is sufficiently
large then
\[
|J_0(\alpha(N)) - I(\alpha(N))| \le \epsilon N.
\]
Then by equation \eqref{eqn:Iak}, for any $\epsilon>0$, if $N$ is
sufficiently large then
\begin{equation}
\label{eqn:J0Est}
J_0(\alpha(N))-J_0(\alpha(0))
\le (2+\epsilon)N.
\end{equation}

On the other hand, by the additivity of $J_0$, we have
\begin{equation}
\label{eqn:J0Sum}
J_0(\alpha(N)) - J_0(\alpha(0)) = \sum_{k=1}^N J_0(C(k)).
\end{equation}
Now it follows from Lemma~\ref{lem:control}(a) that $J_0(C(k))\ge 2$
for all $k$ in the set $A_1\cap A_3$, which has density $1$.  Also
recall from \eqref{eqn:J0Def} and \eqref{eqn:J0ge-1} that $J_0(C(k))$
is always an integer and always at least $-1$.  Combining this with
\eqref{eqn:J0Est} and \eqref{eqn:J0Sum}, we conclude that
$J_0(C(k))=2$ for $k$ in a set of density $1$.
\end{proof}

\subsection{Existence of a genus one Heegaard splitting}
\label{sec:Heegaard}

To complete the proof of Theorem~\ref{thm:lens} (modulo our
temporary simplifying assumption that all Reeb orbits are
nullhomologous), we now show that $Y$ has a genus $1$ Heegaard splitting
such that the Reeb orbits $\gamma_1$ and $\gamma_2$ are the core
circles of the corresponding solid tori.

By Lemma~\ref{lem:prob}, without loss of generality there exists a
cylinder $C\in\mc{M}(\gamma_1^{m_1},\gamma_2^{m_2})$, for some
$m_1,m_2\neq 0$, with ECH index $I(C)=2$.  (It makes no difference in
the argument below if both ends of $C$ are positive.)  Let $\mc{M}_C$
denote the component of $\mc{M}(\gamma_1^{m_1},\gamma_2^{m_2})$
containing $C$.  Let $\pi:\R\times Y\to Y$ denote the projection.  By
Proposition~\ref{prop:foliation}(a), $\pi(C)$ is embedded in $Y$, and
the projections to $Y$ of the cylinders in $\mc{M}_C$ comprise a
foliation of some subset of $Y$.  Since all ECH generators have even
index, it follows from the compactness theorem in \cite[Thm.\
1.8]{pfh2} or \cite[Lem.\ 7.23]{obg1} that the moduli space
$\mc{M}_C/\R$ is compact.  Therefore the projections of the cylinders
in $\mc{M}_C$ foliate all of $Y\setminus (\gamma_1\cup\gamma_2)$, and
$\mc{M}_C/\R\simeq S^1$.

For each $i=1,2$ let $T_i\subset Y$ denote a torus given by the
boundary of a small closed tubular neighborhood $N_i$ of $\gamma_i$.
It follows from the asymptotics for holomorphic curves reviewed in
\S\ref{sec:foliation} that the end at $\gamma_i$ of each holomorphic
curve in $\mc{M}_C$, when projected to $Y$, intersects the torus
$T_i$ transversely in a single circle.  These circles then foliate the
torus $T_i$.

It follows that if $C'\in \mc{M}_C$, then $\pi(C')$ intersects $N_i$
only in a single half-closed cylinder corresponding to the end of $C'$
at $\gamma_i$.  Therefore $\pi(C')$ intersects
$Y\setminus\op{int}(N_1\cup N_2)$ in a closed cylinder.  For each
element of $\mc{M}_C/\R$, we can choose a diffeomorphism of this
closed cylinder with $S^1\times[1,2]$, sending the corresponding
circle in $T_i$ to $S^1\times\{i\}$.  There is no obstruction to
choosing these diffeomorphisms to be smooth functions on
$\mc{M}_C/\R\simeq S^1$, so that they combine to give a diffeomorphism
\[
\varphi: S^1 \times S^1 \times [1,2] \stackrel{\simeq}{\longrightarrow}
Y\setminus \op{int} (N_1\cup N_2)
\]
identifying $S^1\times S^1\times\{i\}$ with $T_i$.  Now
\[
N_1\cup\varphi(S^1\times S^1\times[1,3/2])
\]
and
\[
N_2\cup\varphi(S^1\times S^1\times[3/2,2])
\]
are solid tori in $Y$ which give the desired Heegaard splitting.

\subsection{Removing the simplifying assumption}
\label{sec:RSA}

We now prove Theorem~\ref{thm:lens} without assuming that all Reeb
orbits are nullhomologous.

As in \S\ref{sec:computeI}, let $\gamma_1,\ldots,\gamma_n$ denote the
distinct embedded Reeb orbits.  We know from Lemma~\ref{lem:basic}(b)
that these represent torsion homology classes, so for each
$i=1,\ldots,n$ let $l_i$ denote the smallest positive integer such
that $\gamma_i^{l_i}$ is nullhomologous.  If $m_1,\ldots,m_n$ are
nonnegative integers such that $l_i$ divides $m_i$ for each $i$, then
$\gamma_1^{m_1}\cdots\gamma_n^{m_n}$ is a nullhomologous orbit set.
As such it has a well-defined absolute ECH index.  To compute this,
similarly to \S\ref{sec:computeI}, let $\tau$ be a trivialization of
$\xi$ over the $\gamma_i$'s, let $c_i\eqdef
c_\tau(\gamma_i^{l_i},\phi)/l_i$, let $Q_i\eqdef
Q_\tau(\gamma_i^{l_i},\phi)/l_i^2$, for $i\neq j$ let $Q_{i,j}$ denote
the linking number of $\gamma_i^{l_i}$ and $\gamma_j^{l_j}$ divided by
$l_il_j$, and let $\theta_i$ denote the monodromy angle of $\gamma_i$
with respect to $\tau$.  Then just as in \eqref{eqn:I},
\begin{equation}
\label{eqn:Il}
I(\gamma_1^{m_1}\cdots\gamma_n^{m_n}) = \sum_{i=1}^n(c_im_i +
Q_im_i^2) + \sum_{i\neq j}m_im_jQ_{ij} + \sum_{i=1}^n
\sum_{k=1}^{m_i}(2\floor{k\theta_i}+1).
\end{equation}
In fact, the linear property of $c_\tau$ and the quadratic property of
$Q_\tau$ imply that the index formula \eqref{eqn:Il} is valid for
any nullhomologous ECH generator $\gamma_1^{m_1}\cdots\gamma_n^{m_n}$,
not just one in which each $m_i$ is a multiple of $l_i$.

We can approximate the index formula \eqref{eqn:Il} by
\[
I(\gamma_1^{m_1}\cdots\gamma_n^{m_n}) =
\overline{Q}(m_1,\ldots,m_n) + O(m_1+\cdots+m_n),
\]
where the quadratic form $\overline{Q}$ is defined by
\[
\overline{Q}(m_1,\ldots,m_n) \eqdef
\sum_{i=1}^n(Q_i+\theta_i)m_i^2 + \sum_{i\neq j}m_im_jQ_{ij}.
\]
The same argument as in \S\ref{sec:n=2} then shows that $n=2$.  Also
the same argument as in \S\ref{sec:LBI} shows that Lemma~\ref{lem:LBI}
still holds for nullhomologous ECH generators.  One just needs to
divide each estimate on the number of nullhomologous ECH generators by
the density of the lattice
\[
\left\{(m_1,m_2)\in \Z^2\mid m_1[\gamma_1]+m_2[\gamma_2]=0\in
H_1(Y)\right\}.
\]

Lemma~\ref{lem:prob} then holds by the same argument, using a sequence
of nullhomologus ECH generators provided by
Proposition~\ref{prop:torsion}.  Finally, the argument in
\S\ref{sec:Heegaard} produces a genus $1$ Heegaard splitting just as
before.  Note that $Y\not\simeq S^1\times S^2$ because $[\gamma_1]$
and $[\gamma_2]$ are torsion.  So $Y$ is a lens space, and the proof
of Theorem~\ref{thm:lens} is complete.

\subsection{The theorem on three Reeb orbits}
\label{sec:three}

\begin{proof}[Proof of Theorem~\ref{thm:three}.]
Let $Y$ be a closed oriented connected 3-manifold with a contact form
$\lambda$ such that all Reeb orbits are nondegenerate, and fix a
generic almost complex structure on $\R\times Y$ to define ECH.

First note that if all Reeb orbits are hyperbolic, then there must be
infinitely many distinct embedded Reeb orbits.  Otherwise by
definition the ECH chain complex would have only finitely many
generators, contradicting Proposition~\ref{prop:torsion}.

So by Theorem~\ref{thm:lens}, it is enough to rule out the case where
there are exactly two embedded Reeb orbits, one elliptic and one
hyperbolic.  Suppose that this holds and denote these orbits by $e$
and $h$ respectively.  The ECH generators are now $e^mh^n$ where $m\ge
0$ and $n\in\{0,1\}$.

By Proposition~\ref{prop:torsion}(b), there exists $\Gamma\in H_1(Y)$
such that there are infinitely many ECH generators $\alpha$ with
$[\alpha]=\Gamma$.  It follows that the homology class $[e]\in H_1(Y)$
is torsion.

The homology class $[h]$ is also torsion.  Proof: If
$[h]$ is not torsion then the ECH differential vanishes identically as
in Lemma~\ref{lem:basic}(a).  Since $[e]$ is torsion, it follows that
$ECH_*(Y,\lambda,\Gamma)$ is infinitely generated for both
$\Gamma=0$ and $\Gamma=[h]$.  Since $[h]$ is not torsion, at
least one of these classes $\Gamma$ must have the property that
$c_1(\xi) + 2\op{PD}(\Gamma)\in H^2(Y;\Z)$ is not torsion.  For such a
class $\Gamma$ there are only finitely many possible values of the
grading on $ECH_*(Y,\lambda,\Gamma)$, and now
Proposition~\ref{prop:finite}(a) gives a contradiction.

Since all Reeb orbits represent torsion homology classes, the
cohomology class $c_1(\xi) \in H^2(Y;\Z)$ is then also torsion, as in
Lemma~\ref{lem:basic}(c).

If $[h]\neq 0$, then the same argument as in \S\ref{sec:n=2} gives a
contradiction, by showing that $ECH_*(Y,\lambda,0)$ does not have
enough generators in a given index range in order to be consistent
with Proposition~\ref{prop:torsion}.  If $[h]=0$ then this argument
also works, because $I(e^mh)$ differs from $I(e^m)$ by a linear
function of $m$.
\end{proof}

\section{The Weinstein conjecture for stable Hamiltonian structures}
\label{sec:conclusion}

This section is devoted to the proof of the main Theorem~\ref{thm:main}.

\subsection{Warmup cases}
\label{sec:warmup}

We first prove the theorem in the special cases when $f$ vanishes
either nowhere or everywhere.

If $f$ is nowhere vanishing, then $Y$ is a contact manifold (with the
opposite orientation if $f<0$), so the theorem in this case is just
the Weinstein conjecture which we already know (even if $Y$ is a
$T^2$-bundle over $S^1$).

Now suppose $f=0$, so that the $1$-form $\lambda$ is closed.  Since
$\lambda\wedge\omega > 0$, it follows that $[\lambda]\smile[\omega]\neq
0$ in $H^3(Y;\R)$, so $\lambda$ represents a nonzero cohomology class
in $H^1(Y;\R)$.  We can then add a small closed $1$-form so as to
replace $\lambda$ by a closed $1$-form $\lambda'$ which represents a
real multiple of an integral cohomology class and still satisfies
$\lambda'\wedge\omega>0$ everywhere.  Since $[\lambda']$ is a multiple
of an integral class, there is a fiber bundle $\pi:Y\to S^1$ such that
$\lambda'$ is a multiple of the pullback of the volume form on $S^1$.
Since $\lambda'\wedge\omega>0$, it follows that $\omega$ restricts to
a symplectic form on each fiber of $\pi$.  Since $R$ is in the kernel
of $\omega$, we deduce that $R$ is transverse to the fibers.  If $F$
is a fiber, then the return map of the flow $R$ defines an
orientation-preserving diffeomorphism $\phi:F\to F$, and closed orbits
of $R$ are equivalent to periodic orbits of $\phi$.  By replacing
$\phi$ with an iterate if necessary, we may assume without loss of
generality that $F$ is connected.  The theorem in this case now
follows from part (a) of the following lemma.  Part (b) will be needed
later.

\begin{lemma}
\label{lem:zeta}
Let $F$ be a closed oriented connected surface and let $\phi:F\to F$
be an orientation-preserving diffeomorphism.  Suppose that $\phi$ has
only finitely many irreducible periodic orbits and that all periodic
orbits are nondegenerate and elliptic\footnote{Here ``nondegenerate and
  elliptic'' means that the eigenvalues of the linearized return map
  are not on the real line.}.  Then:
\begin{description}
\item{(a)} If $\phi$ has no periodic orbits, then $F$ is a torus.
\item{(b)} Otherwise $F$ is a sphere, $\phi$ has exactly two
  fixed points, and these are the only irreducible periodic orbits.
\end{description}
\end{lemma}

\begin{proof}
  Let $A$ denote the induced map $\phi_*:H_1(F)\to H_1(F)$.  Since all
  periodic orbits are elliptic, every periodic point of $\phi$ of
  period $p$ counts with weight $+1$ in the Lefschetz fixed point
  formula for $\phi^p$.  We then have the identity
\[
\frac{\det(1-tA)}{(1-t)^2} = \prod_{\gamma}\left(1-t^{p(\gamma)}\right)^{-1}
\]
of formal power series in $\Z[[t]]$, where the product is over
irreducible periodic orbits $\gamma$, and $p(\gamma)$
denotes the period of $\gamma$.  This formula is a special case of the
product formula for the Lefschetz zeta function, see eg \cite{hl},
and it is proved by taking the logarithmic derivative of both sides
and then using the Lefschetz fixed point formula for $\phi$ and its
iterates.  Now we can rewrite the product formula here as
\[
\det(1-tA) \prod_{\gamma}\left(1-t^{p(\gamma)}\right) = (1-t)^2.
\]
Since by hypothesis there are only finitely many factors on the left
hand side, each a polynomial, it follows that the sum of the degrees
of these factors must equal 2, ie
\[
2g(F) + \sum_\gamma p(\gamma) = 2.
\]
The lemma follows immediately.
\end{proof}

To prove Theorem~\ref{thm:main} in the remaining cases, assume that
$R$ has no closed orbit and that the function $f$ is sometimes zero
and sometimes nonzero.  We must show that $Y$ is a $T^2$-bundle over
$S^1$.

\subsection{The region where $f$ is large}
\label{sec:fLarge}

\begin{lemma}
\label{lem:collapse}
  Suppose $\epsilon\ge 0$ is a regular value of $f$.  Then
\[
Y_{\ge \epsilon}\eqdef\{y\in Y
  \mid f(y) \ge \epsilon\}
\]
is diffeomorphic to a disjoint union of
  copies of $T^2\times I$.
\end{lemma}

\begin{proof}
  The idea is to collapse each boundary component of $Y_{\ge
    \epsilon}$ to a circle, so as to obtain a closed contact manifold
  with one embedded elliptic Reeb orbit for each boundary component,
  and then invoke Theorem~\ref{thm:lens}.  We proceed in three steps.

  {\em Step 1.\/} We begin by choosing coordinates near the boundary
  of $Y_{\ge \epsilon}$ in which the stable Hamiltonian structure has
  a nice form.
  
  Fix $\delta>0$ sufficiently small so that every number in the
  interval $[\epsilon,\epsilon+\delta]$ is a regular value of $f$.
  Fix a component $Z$ of $f^{-1}[\epsilon,\epsilon+\delta]$.  For
  $s\in[0,\delta]$, let $\Sigma_s$ denote the component of
  $f^{-1}(\epsilon+s)$ in $Z$.  Each $\Sigma_s$ is a torus, because
  the Reeb vector field $R$ is nonvanishing and tangent to $\Sigma_s$.

  There is a unique vector field $W$ on $Z$ such that $\lambda(W)=0$
  and $\omega(W,\cdot) = df$.  The vector field $W$ is tangent to each
  $\Sigma_s$ and commutes with $R$, and the vectors $R$ and $W$ are
  linearly independent at each point.

  Claim: We can find smooth real-valued functions $\alpha(s), \beta(s),
  \gamma(s), \sigma(s)$ such that for each $s$, the vector fields $\alpha(s)R +
  \beta(s)W$ and $\gamma(s)R + \sigma(s)W$ are linearly independent on
  $\Sigma_s$ and have all orbits closed with period $1$.

  Proof of Claim: Fix $s\in[0,\delta]$.  For $t\in\R$ let
  $\Phi_R^t:\Sigma_s\to\Sigma_s$ and $\Phi_W^t:\Sigma_s\to\Sigma_s$
  denote the time $t$ flows of $R$ and $W$ respectively on $\Sigma_s$.
  Fix a point $p\in\Sigma_s$ and define a map $\phi:\R^2\to \Sigma_s$
  by
\[
\phi(\alpha,\beta)\eqdef
\Phi_R^\alpha\Phi_W^\beta(p).
\]
Since $R$ and $W$ are linearly independent and commute, it follows
that $\phi$ is a covering space and $\phi^{-1}(p)$ is a lattice in
$\R^2$.  Choose $(\alpha(s),\beta(s))$ and $(\gamma(s),\sigma(s))$ to
be a basis for this lattice.  Then $\alpha(s), \beta(s), \gamma(s),
\sigma(s)$ have the required properties for our fixed $s$.  These can
be uniquely extended to smooth functions of $s\in[0,\delta]$ which satisfy the
required properties for all $s$.

It follows from the claim that we can find coordinates
$s\in[0,\delta]$ and $x_1,x_2\in\R/\Z$ on $Z$ in which
\begin{gather*}
f(s,x_1,x_2)=\epsilon+s,\\
R =a_1(s)\frac{\partial}{\partial x_1}
  + a_2(s)\frac{\partial}{\partial x_2}.
\end{gather*}
Since $R$ is assumed to have no closed orbits, the ratio
$a_1(s)/a_2(s)$ is an irrational number which does not depend on $s$.
In fact $a_1(s)$ and $a_2(s)$ do not depend on $s$ either.  To see
this, note that since $\mc{L}_R\lambda = 0$, we have
\[
R\left(\lambda\left(\partial/\partial s\right)\right) 
= - \lambda\left(a_1'(s)\frac{\partial}{\partial x_1} +
a_2'(s)\frac{\partial}{\partial x_2}\right).
\]
Since the right hand side depends only on $s$, and since the function
$\lambda(\partial/\partial s)$ is bounded on $\Sigma_s$, it follows
that the right hand side is zero.  Therefore
$a_1'(s)\frac{\partial}{\partial x_1} +
a_2'(s)\frac{\partial}{\partial x_2}$ is a multiple of $W$, and this
multiple must be zero because $a_1(s)/a_2(s)$ is constant.

Since $\omega(R,\cdot)=0$, we can write
\[
\omega = c(s)ds\wedge(a_2 dx_1 - a_1 dx_2).
\]
Note that the function $c$ depends only on $s$, because
$\mc{L}_R\omega=0$ implies $R(c)=0$, and $R$ is ergodic on each
$\Sigma_s$.  Also $c$ is nowhere vanishing, and we can choose our
coordinates so that $c$ is always positive.  Likewise, since
$\lambda(R)=1$, we can write
\[
\lambda = b(s)(a_2dx_1-a_1dx_2) +
\sigma(s)ds+(a_1^2+a_2^2)^{-1}(a_1dx_1+a_2dx_2),
\]
where $b$ and $\sigma$ depend only on $s$ because $\mc{L}_R\lambda=0$.
Since $d\lambda = f\omega$, we have
\begin{equation}
\label{eqn:bprime}
b'(s)=(\epsilon + s)c(s).
\end{equation}
Also, changing $\lambda$ by adding the differential of some function
of $s$ does not change its salient properties, so we may assume that
$\sigma(s)=0$.

{\em Step 2.\/} We now collapse the boundary component $\Sigma_0$ to
a circle.

To prepare for this, define a function $\rho:[0,\delta]\to\R$ by
\[
\rho(s)^2 \eqdef 2(b(s)-b(0)).
\]
Since $\epsilon\ge 0$ and $c(s)$ is always positive, it follows from
equation \eqref{eqn:bprime} that $\rho$ is strictly increasing and smooth
on $(0,\delta]$.  In terms of this function we can write
\[
\lambda = \frac{1}{2}\rho^2(a_2dx_1 - a_1dx_2) + \alpha_1dx_1 +
\alpha_2dx_2
\]
where $\alpha_1$ and $\alpha_2$ are constants satisfying $a_1\alpha_1 +
a_2\alpha_2 = 1$.

We next modify $\lambda$ on $Y_{\ge\epsilon}$, without changing its
salient properties, to arrange that $\alpha_1/\alpha_2$ is rational.
Suppose that $\alpha_1/\alpha_2$ is irrational.  The restriction map
$H^1(Y_{\ge\epsilon};\Z)\to H^1(\Sigma_0;\Z)$ is nonzero, and so there
exists a closed $1$-form $\lambda_1$ on $Y_{\ge\epsilon}$ such that
$\lambda_1|_{Z}=\beta_1dx_1+\beta_2dx_2$ for some relatively prime
integers $\beta_1$ and $\beta_2$.  Now consider replacing $\lambda$ by
$\lambda + \tau\lambda_1$ where $\tau$ is a small constant.  If
$\tau>0$ is sufficiently small then we still have $(\lambda +
\tau\lambda_1)\wedge\omega>0$.  And since $\alpha_1/\alpha_2$ is
irrational and $\beta_1/\beta_2$ is rational, it follows that
$(\alpha_1+\tau\beta_1)/(\alpha_2+\tau\beta_2)$ is rational for a
dense set of $\tau$.  Replacing $\lambda$ by $\lambda+\tau\lambda_1$
multiplies the Reeb vector field by a positive function; in
particular there are still no Reeb orbits.

With $\alpha_1/\alpha_2$ arranged to be rational, by an $SL_2\Z$
coordinate change we can further assume that $\alpha_1=0$.  Now let
$Y_{\ge\epsilon}'$ be obtained from $Y_{\ge\epsilon}$ by declaring two
points in $\Sigma_0$ to be equivalent whenever they have the same
$x_1$ coordinate.  The subset $Z$ of $Y_{\ge\epsilon}$ gets collapsed
to a subset $Z'$ of $Y_{\ge\epsilon}'$ which is a disk cross $S^1$.
Define the smooth structure on $Z'$ so that $(\rho,x_1)$ are polar
coordinates on the disk and $x_2$ is the $S^1$ coordinate.  Then
$\lambda$ is a smooth contact form on $Z'$, because in the above
coordinates,
\[
\lambda  = \frac{1}{2}\rho^2(a_2dx_1 - a_1dx_2) + \alpha_2dx_2
\]
is smooth, and
\[
\lambda\wedge d\lambda = \rho\, d\rho\, dx_1\, dx_2
\]
is nonvanishing.  The boundary torus $\Sigma_0$ gets collapsed to a
new Reeb orbit $\Sigma_0'$ in $Y_{\ge\epsilon}'$.  This Reeb orbit is
elliptic with irrational monodromy angle, because $R$ has irrational
slope on the tori where $\rho$ is constant.  The Reeb vector field on
$Y_{\ge\epsilon}'\setminus\Sigma_0' =
Y_{\ge\epsilon}\setminus\Sigma_0$ is unaffected.

{\em Step 3.\/} We now complete the proof of Lemma~\ref{lem:collapse}.

Repeating the above process finitely many times, we collapse all
boundary components of $Y_{\ge\epsilon}$ to circles, to obtain a
closed contact manifold with one embedded elliptic Reeb orbit for each
boundary component, all of whose iterates are nondegenerate, and no
other embedded Reeb orbits.  Let $X$ be a component of
$Y_{\ge\epsilon}$.  It then follows from
Theorem~\ref{thm:lens} that there are exactly two such orbits which
comprise the core circles of a genus $1$ Heegaard splitting of the
collapsed $X$.  Uncollapsing these circles back to
boundary tori, we conclude that $X$ is diffeomorphic to
$T^2\times I$.
\end{proof}

The same argument with some different signs shows that if $\epsilon\le
0$ is a regular value of $f$, then
\[
Y_{\le\epsilon} \eqdef \{y\in Y\mid f(y)\le \epsilon\}
\]
is likewise diffeomorphic to a disjoint union of copies of $T^2\times
I$.  In particular, if $0$ is a regular value of $f$, then it follows
that $Y$ is a union of copies of $T^2\times I$ glued together along
their boundaries, so Theorem~\ref{thm:main} is proved in this case.

\subsection{The region where $f$ is small}
\label{sec:fSmall}

We now begin the proof of Theorem~\ref{thm:main} when
$0$ is not necessarily a regular value of $f$.  Choose a small
$\epsilon>0$ such that both $\epsilon$ and $-\epsilon$ are regular
values of $f$.
Define
\[
Y_{\epsilon} \eqdef \{y\in Y \mid |f(y)| \le \epsilon\}.
\]
To prove Theorem~\ref{thm:main}, we will show that if $\epsilon$ as
above is sufficiently small, then $Y_{\epsilon}$ is a disjoint union
of copies of $T^2\times I$.  The strategy for doing so is to perturb
$\lambda|_{Y_\epsilon}$ to a closed form which still has positive
wedge product with $\omega$, deduce that $Y_\epsilon$ fibers over
$S^1$ such that the Reeb vector field is transverse to the fibers, and
then apply Lemma~\ref{lem:zeta}(b) to show that the fibers are
disjoint unions of annuli.

To start, fix a metric on $Y$.  Also fix a smooth function
$\chi:[0,1]\to[0,1]$ such that $\chi(t)=1$ for $t<1/3$ and $\chi(t)=0$
for $t>2/3$.

Fix $\epsilon$ as above.  Choose $\delta\in(0,\epsilon)$ sufficiently
small so that all numbers in the intervals
$[\epsilon,\epsilon+\delta]$ and $[-\epsilon-\delta,-\epsilon]$ are
regular values of $f$.  Define a $1$-form $\lambda_*$ on
$Y_{\epsilon+\delta}$ as follows.  On $Y_\epsilon$ define
$\lambda_*\eqdef\lambda$.  We now define $\lambda_*$ on
$f^{-1}[\epsilon,\epsilon+\delta]$.  The construction on
$f^{-1}[-\epsilon-\delta,-\epsilon]$ is analogous and will be omitted.

Fix a component $Z$ of $f^{-1}[\epsilon,\epsilon+\delta]$.  Recall
from \S\ref{sec:fLarge} that we can choose coordinates
$s\in[0,\delta]$ and $x_1,x_2\in\R/\Z$ on $Z$, and modify $\lambda$
without changing its salient properties, so that in these coordinates
\begin{align*}
\lambda &= b(s)(a_2dx_1-a_1dx_2) +
\frac{a_1dx_1+a_2dx_2}{a_1^2+a_2^2},\\
\omega &= c(s)ds\wedge(a_2dx_1-a_1dx_2),
\end{align*}
where $c(s)>0$.  We now define $\lambda_*$ on $Z$ by
\[
\lambda_* \eqdef \left[b(0)+\chi(\delta^{-1}s)(b(s)-b(0)\right]
(a_2dx_1-a_1dx_2) +
\frac{a_1dx_1+a_2dx_2}{a_1^2+a_2^2}.
\]
On $Y_{\epsilon+\delta}$ we then have
\[
d\lambda_* = g\omega
\]
where $g$ is a smooth function which agrees with $f$ on $Y_\epsilon$
and extends by zero to a smooth function defined on all of $Y$.

\subsection{An upper bound on $g$}

We now show that there is an $\epsilon$-independent constant $c_0$
such that if $\delta$ is chosen sufficiently small in the above
construction then
\begin{equation}
\label{eqn:gBound}
|g|\le c_0\epsilon.
\end{equation}

We just need to check this on a region $Z$ as above.  It follows from
the above equations that on $Z$ we have
\begin{equation}
\label{eqn:gZ}
g(s) = c(s)^{-1}\left[\delta^{-1}\chi'(\delta^{-1}s)(b(s)-b(0)) +
  \chi(\delta^{-1}s)b'(s)\right].
\end{equation}
Since $b'(s)=(\epsilon+s)c(s)$ and we have chosen $\delta<\epsilon$,
it follows that
\[
c(s)^{-1}\chi(\delta^{-1}s)b'(s) =
\chi(\delta^{-1}s)(\epsilon+s) < 2\epsilon.
\]
If we further choose $\delta$ sufficiently small so that $c(t)\le
2c(s)$ for all $s,t\in[0,\delta]$, then we can estimate
\[
\begin{split}
b(s)-b(0) &= \int_{0}^s (\epsilon + t)c(t)dt\\
&\le 2c(s)(\epsilon s + s^2/2)\\
&\le 3 c(s)\epsilon\delta.
\end{split}
\]
Putting the above two inequalities into \eqref{eqn:gZ} proves
\eqref{eqn:gBound}.

\subsection{The cohomology class of $g\omega$}

We now study the cohomology class of the closed
$2$-form $g\omega$ on $Y$.

Let $U\subset Y$ denote the open set where $f\neq 0$.  Consider the
relative homology long exact sequence
\[
\cdots \longrightarrow H_1(U) \longrightarrow H_1(Y) \longrightarrow
H_1(Y,U) \longrightarrow \cdots.
\]
Here and below all homology and cohomology is with real coefficients.
Fix embedded oriented curves $\gamma_1,\ldots,\gamma_n\subset U$ that
represent a basis for the kernel of the map $H_1(Y)\to H_1(Y,U)$.
Since the curves $\gamma_i$ all have positive distance from the
compact set where $f=0$, there exists $d>0$ such that if $\epsilon$ is
sufficiently small then no point in $Y_{\epsilon}$ is within distance
$d$ of any point in any of the curves $\gamma_i$.  For each
$i=1,\ldots,n$, fix an $\epsilon$-independent closed $2$-form $\mu_i$
which represents the Poincar\'{e} dual of $\gamma_i$ and is supported
within distance $d/2$ of $\gamma_i$.

Claim: there exists an $\epsilon$-independent constant $c_0$ such that
if $\epsilon$ is chosen sufficiently small in the construction in
\S\ref{sec:fSmall}, then there are unique real numbers
$q_1,\ldots,q_n$ such that
\begin{equation}
\label{eqn:qi}
[g\omega] = \sum_{i=1}^nq_i[\mu_i] \in H^2(Y),
\end{equation}
and these satisfy
\begin{equation}
\label{eqn:qiBound}
|q_i|\le c_0\epsilon.
\end{equation}
Note that $c_0$ here is different than in \eqref{eqn:gBound}.

Proof of claim: Choose $\epsilon$ sufficiently small so that
$Y_\epsilon$ does not intersect the support of the forms $\mu_i$. Then the
Poincar\'{e} duality isomorphism $H_1(Y)\stackrel{\simeq}{\to} H^2(Y)$
restricts to an injection
\begin{equation}
\label{eqn:PDR}
\Ker(H_1(Y)\to H_1(Y,U)) \longrightarrow \Ker(H^2(Y)\to
H^2(Y_\epsilon)).
\end{equation}

Since $q\omega$ restricts to an exact form on $Y_\epsilon$, to prove
that the numbers $q_i$ exist and are unique it is enough to show that
the map \eqref{eqn:PDR} is surjective.  To prove this surjectivity,
note that any element of $\Ker(H^2(Y)\to H^2(Y_\epsilon))$ can be
represented by a closed $2$-form $\eta$ with support on $Y\setminus
Y_\epsilon$.  Here we are assuming as usual that $\epsilon$ and
$-\epsilon$ are regular values of $f$.  Now the Poincare dual of
$[\eta]\in H^2_c(Y\setminus Y_\epsilon)$ is a homology class
$\alpha\in H_1(Y\setminus Y_\epsilon)$ with $\alpha\cdot S =
\int_S\eta$ for all $S\in H_2(Y\setminus Y_\epsilon, \partial
Y_\epsilon)$.  Letting $\imath:Y\setminus Y_\epsilon\to Y$ denote the
inclusion, we then have $\imath_*\alpha \cdot S = \int_S\eta$ for all
$S\in H_2(Y)$.  So $\imath_*\alpha\in H_1(Y)$ is the Poincar\'{e} dual
of $[\eta]\in H^2(Y)$, and $\imath_*\alpha$ maps to $0$ in $H_1(Y,U)$
since $Y\setminus Y_\epsilon\subset U$.

To prove that $|q_i|\le c_0\epsilon$, for each $i=1,\ldots,n$ we can fix an
$\epsilon$-independent closed oriented embedded surface $\Sigma_i$ in
$Y$ which has intersection number $\delta_{ij}$ with $\gamma_j$.  Then
observe that
\[
q_i = \int_{\Sigma_i}\sum_{j=1}^n q_j\mu_j = \int_{\Sigma_i}g\omega
\]
and use \eqref{eqn:gBound}.

\subsection{The region where $f$ is small fibers over $S^1$}
\label{sec:fibration}

We now show that if $\epsilon$ is sufficiently small
then $Y_\epsilon$ fibers over $S^1$ so that the Reeb vector field $R$
is transverse to the fibers.

Let $\epsilon>0$ be small enough so that $Y_\epsilon$ does not
intersect the support of the forms $\mu_i$.  Let $q_i$ be the numbers
satisfying \eqref{eqn:qi}.  Since the $2$-form $g\omega -
\sum_{i=1}^nq_i\mu_i$ is exact on $Y$, there exists a unique $1$-form
$\nu$ such that $d\nu = g\omega - \sum_{i=1}^n q_i\mu_i$ and
$d^*\nu=0$ and $\nu$ is $L^2$-orthogonal to the space of harmonic
$1$-forms on $Y$.  It then follows from estimates on the Green's
function for $d+d^*$, namely the fact that the singularity of the
Green's function behaves as $\op{dist}(y_1,y_2)^{-2}$ for any points
$y_1\neq y_2$ in $Y$, that there is an $\epsilon$-independent constant
$c_0$ such that
\[
|v|\le c_0\sup_Y\left|g\omega - \sum_{i=1}^nq_i\mu_i\right|.
\]
Together with \eqref{eqn:gBound} and \eqref{eqn:qiBound}, this implies
that there is an $\epsilon$-independent constant $c_0$ such that
\begin{equation}
\label{eqn:nuBound}
|\nu| \le c_0\epsilon.
\end{equation}

Now the $1$-form $\lambda - \nu$ restricts to a closed $1$-form on
$Y_\epsilon$.  By the estimate \eqref{eqn:nuBound}, if $\epsilon$ is
sufficiently small then $(\lambda-\nu)\wedge\omega>0$ everywhere.  We
can perturb $\lambda-\nu$ to a closed $1$-form $\lambda'$ on
$Y_\epsilon$ which represents a real multiple of an integral
cohomology class in $H^1(Y_\epsilon;\Z)$ and still satisfies
$\lambda'\wedge\omega>0$.  It follows as in \S\ref{sec:warmup} that
$Y_\epsilon$ fibers over $S^1$ with the Reeb vector field $R$
transverse to the fibers.

\subsection{The fibers are disjoint unions of annuli}

To complete the proof of Theorem~\ref{thm:main}, we now show that
$Y_\epsilon$ is diffeomorphic to a disjoint union of copies of
$T^2\times I$.

Let $X$ be a component of $Y_\epsilon$.  Let $F$ be a fiber of the
fibration $X\to S^1$ constructed in \S\ref{sec:fibration}, and let
$\phi:F\to F$ denote the return map of the flow $R$.  Closed orbits of
$R$ in $X$ are equivalent to periodic orbits of $\phi$, so by
assumption $\phi$ has no periodic orbits.  By replacing $\phi$ by an
iterate if necessary, we may assume without loss of generality that
$F$ is connected.

Assume as usual that $\epsilon$ and $-\epsilon$ are regular values of
$f$, so that the boundary of $X$ is a disjoint union of tori.  Each
boundary circle of $F$ lies in a boundary torus of $X$.  Recall from
\S\ref{sec:fLarge} that on each boundary torus of $X$, in suitable
coordinates the Reeb vector field $R$ is a constant vector field with
irrational slope which does not change if we perturb $\epsilon$.

Let $F'$ be the closed surface obtained from $F$ by identifying two
points whenever they are in the same boundary circle.  Then $\phi$
descends to a diffeomorphism $\phi':F'\to F'$ with one irreducible
periodic orbit $\gamma$ for each boundary torus of $X$, and no other
irreducible periodic orbits.  The period of $\gamma$ equals the number
of boundary circles of $F'$ on the corresponding boundary torus of
$X$.  It follows from the above description of $R$ near the boundary
of $X$ that the orbit $\gamma$ and all of its iterates are
nondegenerate and elliptic.

We now invoke Lemma~\ref{lem:zeta}(b) to conclude $F'$ is a sphere and
$\phi'$ has two fixed points and no other irreducible periodic orbits.
It follows immediately that $F$ is an annulus, so $X$ is diffeomorphic
to $T^2\times I$ as desired.  This completes the proof of
Theorem~\ref{thm:main}.

\end{document}